\documentclass{article}%
\usepackage{amsmath}
\usepackage{amssymb}
\usepackage{amsfonts}
\usepackage{graphicx}%
\newtheorem{Th}{Theorem}[section]
\newtheorem{Prop}[Th]{Proposition}
\newtheorem{Lemma}[Th]{Lemma}
\newtheorem{Coro}[Th]{Corollary}

\newtheorem{Rem}[Th]{Remark}

\begin{document}

\date{}
\title{Solutions of mKdV in classes of functions unbounded at infinity}
\author{T. Kappeler\thanks{Supported in part by the Swiss National Science Foundation,
and the programme SPECT, and the European Community through the FP6 Marie
Curie RTN ENIGMA (MRTN-CT-2004-5652).}, P. Perry\thanks{Partially supported by
NSF-grant DMS-0408419}, M. Shubin\thanks{Partially supported by NSF-grant
DMS-0600196}, and P. Topalov}
\maketitle

\begin{abstract}
In 1974 P. Lax introduced an algebro-analytic mechanism similar to the
Lax L-A pair. Using it we prove global existence and uniqueness for solutions
of the initial value problem for mKdV in classes of smooth functions which can be
unbounded at infinity, and may even include functions which tend to infinity with respect
to the space variable. Moreover, we establish the invariance of the spectrum
and the unitary type of the Schr{\"o}dinger operator under the KdV flow and the invariance of the
spectrum and the unitary type of the impedance operator under the mKdV flow
for potentials in these classes.

\vspace{0.3cm}

\noindent{\em Mathematics Subject Classification 2000:} 34A12, 35053, 37K40

\vspace{0.1cm}

\noindent{\em Keywords:} KdV, modified KdV, spectra of Schr{\"o}dinger operators

\end{abstract}

\section{Introduction}

The purpose of this work is to solve the \emph{modified Korteweg - de Vries
equation} (mKdV) on the line
\begin{align}
r_{t}-6r^{2}r_{x}+r_{xxx}=0\label{mKdV1}\\
r|_{t=0}=r_{0} \label{mKdV2}%
\end{align}
in various classes of smooth functions (possibly) unbounded at $+\infty$
and/or $-\infty$. Equation (\ref{mKdV1}) is closely related to the celebrated
Korteweg - de Vries equation (KdV),
\begin{equation}
\label{KdV}q_{t}-6qq_{x}+q_{xxx}=0
\end{equation}
and is a model equation for wave propagation.

Let $I=(a,b)\subseteq\mathbb{R}$ with $-\infty\le a<b\le\infty$. For any given
$\beta\in\mathbb{R}$ denote by $\mathcal{S}_{\beta}(I\times\mathbb{R})$ the
linear space of $C^{\infty}(I\times\mathbb{R})$ functions having asymptotic
expansions at $+\infty$ and $-\infty$ (cf. \cite{BS1})
\begin{equation}
\label{at+infty}r(t,x)\sim\sum_{k=0}^{\infty}a_{k}^{+}(t)x^{\beta_{k}%
}\;\;\mbox{as}\;\;x\to\infty
\end{equation}
and
\begin{equation}
\label{at-infty}r(t,x)\sim\sum_{k=0}^{\infty}a_{k}^{-}(t)(-x)^{\beta_{k}%
}\;\;\mbox{as}\;\;x\to-\infty
\end{equation}
where $a_{k}^{\pm}\in C^{\infty}(I)$ and $\beta=\beta_{0}>\beta_{1}>...$ with
$\lim\limits_{k\to\infty}\beta_{k}=-\infty$. By definition, the relations
(\ref{at+infty}) and (\ref{at-infty}) mean that for any compact interval
$J\subseteq I$ and any $N\ge0$, $i,j\ge0$, there exists a constant
$C_{J,N,i,j}>0$ such that for any $\pm x\ge1$ and $t\in J$
\begin{equation}
\label{e:asymptotics}\Big|\;\partial_{t}^{i}\partial_{x}^{j}\Big(r(t,x)-
\sum_{k=0}^{N} a_{k}^{\pm}(t)(\pm x)^{\beta_{k}}\Big)\Big|\le C_{J,N,i,j}%
|x|^{\beta_{N+1}-j}.
\end{equation}
For an arbitrarily chosen formal series $\sum_{k=0}^{\infty}a_{k}^{\pm}(t)(\pm
x)^{\beta_{k}}$, referred to as a symbol in the theory of pseudodifferential
operators, there exists a function $r\in C^{\infty}(I\times\mathbb{R})$
satisfying (\ref{at+infty}) and (\ref{at-infty}) (see for example
\cite[Proposition 3.5]{Shubin}). Analogously one defines the linear space
$\mathcal{S}_{\beta}(\mathbb{R})$ as the space of functions $r\in C^{\infty
}(\mathbb{R})$ having asymptotic expansions $r(x)\sim\sum_{k=0}^{\infty}%
a_{k}^{\pm}(\pm x)^{\beta_{k}}\;\;\mbox{as}\;\;x\to\pm\infty$ where
$a_{k}^{\pm}$ are given constants, $\beta=\beta_{0}>\beta_{1}>...$ and
$\lim\limits_{k\to\infty}\beta_{k}=-\infty$.

In this paper we first prove the following results about the initial value
problem \eqref{mKdV1}-\eqref{mKdV2}.

\begin{Th}
\hspace{-2mm}\textbf{.}\label{Th:mKdV-S} For any $\beta<1/2$ and for any
initial data $r_{0}\in\mathcal{S}_{\beta}(\mathbb{R})$ there exists a solution
$r\in\mathcal{S}_{\beta}(\mathbb{R}\times\mathbb{R})$ of the initial value
problem (\ref{mKdV1})-(\ref{mKdV2}).
The solution $r$ is unique in the class of solutions of (\ref{mKdV1}%
)-(\ref{mKdV2}) in $\mathcal{S}_{\beta}(\mathbb{R}\times\mathbb{R})$.
Moreover, the coefficients $a_{0}^{\pm}(t)$ in the asymptotic expansion of the
solution $r(t,x)$ are independent of $t$ and are equal to the coefficients
$a_{0}^{\pm}$ in the asymptotic expansion of the initial data $r_{0}$.
\end{Th}

Note that the solution $r$ is global in time which will be also the case for
all results formulated below.

By the same method of proof we obtain similar results for the larger spaces of
functions $\mathcal{O}_{\beta}(I\times\mathbb{R})$ and $o_{\beta}%
(I\times\mathbb{R})$ which are (possibly) unbounded at infinity.

Let $I=(a,b)\subseteq\mathbb{R}$ with $-\infty\le a<b\le\infty$. For any given
$\beta\in\mathbb{R}$ denote by $\mathcal{O}_{\beta}(I\times\mathbb{R})$ the
linear space of functions $r(t,x)$ in $C^{\infty}(I\times\mathbb{R})$ such
that for any compact interval $J\subseteq I$ and any $k,l\ge0$ there exists a
constant $C_{J,k,l}>0$ such that for any $|x|\ge1$ and any $t\in J$
\[
|\partial_{t}^{k}\partial_{x}^{l}r(t,x)|\le C_{J,k,l}|x|^{\beta-l}.
\]
Analogously one defines the linear space $\mathcal{O}_{\beta}(\mathbb{R})$ as
the space of functions $r(x)$ in $C^{\infty}(\mathbb{R})$ such that for any
$l\ge0$ there exists $C_{l}>0$ such that for any $|x|\ge1$, $|\partial_{x}%
^{l}r(x)|\le C_{l}|x|^{\beta-l}$.

We will also consider the following spaces. For any given $\beta\in\mathbb{R}$
denote by ${o}_{\beta}(I\times\mathbb{R})$ the linear space of functions
$r(t,x)$ in $C^{\infty}(I\times\mathbb{R})$ such that for any compact interval
$J\subseteq I$ and any $k,l\ge0$
\[
\partial_{t}^{k}\partial_{x}^{l}r(t,x)=o(|x|^{\beta-l})
\]
uniformly in $t\in J$. In the same way as above one defines the space
${o}_{\beta}(\mathbb{R})$. Clearly the following inclusions hold:
\[
\mathcal{S}_{\beta}(I\times\mathbb{R})\subseteq\mathcal{O}_{\beta}%
(I\times\mathbb{R}),\;\;\; {o}_{\beta}(I\times\mathbb{R})\subseteq
\mathcal{O}_{\beta}(I\times\mathbb{R}).
\]

\begin{Th}
\hspace{-2mm}\textbf{.}\label{Th:mKdV-O} For any $\beta<1/2$ and for any
initial data $r_{0}\in\mathcal{O}_{\beta}(\mathbb{R})$ there exists a global
in time solution $r\in\mathcal{O}_{\beta}(\mathbb{R}\times\mathbb{R})$ a
solution $r\in\mathcal{O}_{\beta}(\mathbb{R}\times\mathbb{R})$ of the initial
value problem (\ref{mKdV1})-(\ref{mKdV2}).
The solution $r$ is unique in the class of solutions of (\ref{mKdV1}%
)-(\ref{mKdV2}) in $\mathcal{O}_{\beta}(\mathbb{R}\times\mathbb{R})$.
\end{Th}

\begin{Th}
\hspace{-2mm}\textbf{.}\label{Th:mKdV-o} For any $\beta\le1/2$ and for any
initial data $r_{0}\in{o}_{\beta}(\mathbb{R})$ there exists a solution
$r\in{o}_{\beta}(\mathbb{R}\times\mathbb{R})$ of the initial value problem
(\ref{mKdV1})-(\ref{mKdV2}). The solution $r$ is unique in the class of
solutions of (\ref{mKdV1})-(\ref{mKdV2}) in ${o}_{\beta}(\mathbb{R}%
\times\mathbb{R})$.
\end{Th}

\begin{Rem}
\hspace{-2mm}\textbf{.}\label{optimal} Note that for $r_{0} \in\mathcal{S}%
_{\beta}(\mathbb{R})$ with $\beta= \beta_{0} > 1/2,$ with an asymptotic
expansion of the form
\[
r_{0}(x)\sim\sum_{k=0}^{\infty}a_{k}^{+}x^{\beta_{k}}\;\;\mbox{as}\;\;x\to
+\infty
\]
with $a_{0}^{+} \ne0$, no formal solution and therefore no solution of mKdV in
$\mathcal{S}_{\beta}(\mathbb{R} \times\mathbb{R})$ exists.
\end{Rem}

\begin{Rem}
\hspace{-2mm}\textbf{.} In fact, the uniqueness in all Theorems
\ref{Th:mKdV-S}, \ref{Th:mKdV-O}, \ref{Th:mKdV-o} holds if we only require
$r\in{o}_{1/2}(\mathbb{R}\times\mathbb{R})$ which is the largest class where
the existence is claimed in these theorems. So if we only require $r\in
{o}_{1/2}(\mathbb{R}\times\mathbb{R})$ and take the initial condition $r_{0}$
in $\mathcal{S}_{\beta}(\mathbb{R})$, $\mathcal{O}_{\beta}(\mathbb{R})$
$(\beta<1/2)$ or ${o}_{\beta}(\mathbb{R})$ $(\beta\le1/2)$, then we will
automatically have $r\in\mathcal{S}_{\beta}(\mathbb{R}\times\mathbb{R})$,
$\mathcal{O}_{\beta}(\mathbb{R}\times\mathbb{R})$ or ${o}_{\beta}%
(\mathbb{R}\times\mathbb{R})$ respectively.
\end{Rem}

Results for KdV similar to the ones stated for mKdV in Theorems
\ref{Th:mKdV-S}, \ref{Th:mKdV-O}, and \ref{Th:mKdV-o} have been obtained in a
series of papers \cite{B1,BS1,BS2,BS3} -- see Appendix B where, for the
convenience of the reader, we give a short summary of these results. In fact,
we construct our solutions of mKdV with the properties stated in the above
theorems by applying to the solutions of \cite{B1,BS1,BS2,BS3} an inverse of
the Miura map. Recall that the Miura map $r\mapsto B(r):=r_{x}+r^{2}$, first
introduced in \cite{Miura}, maps smooth solutions of mKdV to smooth solutions
of KdV. However, the Miura map is usually neither 1-1 nor onto. This is, for
example, the case when $B$ is considered as a map $H_{loc}^{\beta}%
(\mathbb{R})\rightarrow H_{loc}^{\beta-1}(\mathbb{R})$ with $\beta\geq0$
\cite{KPST}. In this case, the preimage of an element in $H_{loc}^{\beta
-1}(\mathbb{R})$ is either the empty set, a point or a set homeomorphic to an
interval. To describe the preimage $B^{-1}\{B(r)\}$ of $q=B(r)$, note that the
positive function $\psi(x)=e^{\int_{0}^{x}r(s)ds}$ satisfies
\begin{equation}
-\psi_{xx}+(r_{x}+r^{2})\psi=0 \label{schroed}%
\end{equation}
and is related to $r$ by $r=\psi_{x}/\psi$. It has been shown in \cite{KPST}
that for $r\in H_{loc}^{\beta}(\mathbb{R})$ given with $\beta\geq0$, any
function in the preimage $B^{-1}\{B(r)\}$ arises in this way, i.e. for any
$r\in H_{loc}^{\beta}(\mathbb{R}),$
\[
B^{-1}\{B(r)\}=\{\psi_{x}/\psi\mid\psi\in H_{loc}^{\beta}(\mathbb{R})\text{
positive, satisfying }(\ref{schroed})\}.
\]
Given initial data $r_{0}$ in the class of functions considered in the
theorems above, $q_{0}=B(r_{0})$ has the growth condition at infinity required
by the theorems in \cite{B1,BS1,BS2,BS3} to conclude that there exists a
unique solution $q(t,x)$ of KdV in the corresponding class with $q(0,\cdot
)=q_{0}$. We then consider the linear evolution equation, introduced by Lax
\cite{Lax:1974} (see also Marchenko \cite{Mar}),
\begin{align}
\psi_{t}(t,x)  &  =Q(t)\psi(t,x)\label{evol1}\\
\psi(0,x)  &  =e^{\int_{0}^{x}r_{0}(s)ds} \label{evol2}%
\end{align}
where $Q(t)$ is the first-order differential operator,
\begin{equation}
Q(t):=2q(t,x)\partial_{x}-q_{x}(t,x). \label{opQ}%
\end{equation}
and prove that there exists a unique, globally (in time) defined solution
$\psi(t,x)$, satisfying $\psi(t,x)>0$ for any $x\in\mathbb{R},$ $t\in
\mathbb{R}$ and
\begin{equation}
-\psi_{xx}(t,x)+q(t,x)\psi(t,x)=0. \label{schroed-t}%
\end{equation}
The latter identity follows from the commutator relation
\begin{equation}
{\dot{L}}=[Q,L]+4q_{x}L \label{ldot}%
\end{equation}
where
\begin{equation}
L(t):=-\partial_{x}^{2}+q(t,x). \label{opL}%
\end{equation}
and ${\dot{L}}=q_{t}.$ The function
\begin{equation}
r(t,x):=\psi_{x}(t,x)/\psi(t,x) \label{representation}%
\end{equation}
is then the unique solution of mKdV with $r(0,\cdot)=r_{0}$ in a class of
functions in $C^{\infty}(\mathbb{R}\times\mathbb{R})$ satisfying appropriate
growth conditions. It has the claimed properties in each of the settings of
Theorem \ref{Th:mKdV-S}, \ref{Th:mKdV-O}, and \ref{Th:mKdV-o}. The pair of
operators $(Q,L)$ can be obtained from the classical Lax L-A pair for KdV by
division--see section 2 below. We refer to it as a Q-L pair. Such a pair
allows us to construct an inverse of the Miura map and, in this way, deduce
existence and uniqueness of solutions for (\ref{mKdV1})-(\ref{mKdV2}) from the
corresponding results for KdV.

We also establish the invariance of the spectrum of the Schr{\"o}dinger
operator under the KdV flow and the invariance of the spectrum of the
impedance operator under the mKdV flow. Consider the Schr{\"o}dinger operator
$L(t)=-\frac{d^{2}}{dx^{2}}+q(t,x)$ where $q(t,x)$ is a solution of the KdV
equation in $\mathcal{O}_{\beta}(I\times\mathbb{R})$ with $\beta\le1$ and
$I=(a,b)$, $-\infty\le a<b\le+\infty$. (By Sears theorem (cf. \cite[Chapter
II]{BeS}), for any given $t\in(a,b)$ the operator $-\frac{d^{2}}{dx^{2}%
}+q(t,x)$ with domain $C^{\infty}_{0}(\mathbb{R})$ is essentially
self-adjoint. Denote by $L(t)$ its closure.) We will prove the following

\begin{Th}
\hspace{-2mm}\textbf{.}\label{Prop:spectral_invariance} Let $q\in
\mathcal{O}_{\beta}(I\times\mathbb{R})$ with $\beta\le1$ be a solution of the
KdV equation.

\begin{itemize}
\item[(i)] Then for any $t,t^{\prime}\in I$
\[
\mathop{\tt spec} L(t)= \mathop{\tt spec} L(t^{\prime})\,.
\]
Moreover, the point spectra of the operators $L(t)$ and $L(t^{\prime})$ (i.e.
the sets of eigenvalues corresponding to $L^{2}$-eigen\-functions), coincide
and have the same multiplicities.

\item[(ii)] If in addition $\beta<1$ or $q\in o_{\beta}(I\times\mathbb{R})$
with $\beta\le1$, then for any $t, t^{\prime}\in I$ the operators $L(t)$ and
$L(t^{\prime})$ are unitarily equivalent. It means that for any $t,t^{\prime
}\in I$ there exists a unitary operator $\Psi(t,t^{\prime}) : L^{2}%
(\mathbb{R})\to L^{2}(\mathbb{R})$ such that $\Psi(t^{\prime}%
,t)L(t)=L(t^{\prime})\Psi(t^{\prime},t)$.
\end{itemize}
\end{Th}

\begin{Rem}
\hspace{-2mm}\textbf{.} Our method suggests that the operators $L(t)$ and
$L(t^{\prime})$ are unitarily equivalent also in the case when $q\in
\mathcal{O}_{\beta}(I\times\mathbb{R})$ with $\beta\le1$. However, the proof
of this statement will require an improvement of Theorem 2 in \cite{BS2}.
\end{Rem}

In order to prove the first statement of Theorem
\ref{Prop:spectral_invariance} we again use the Q-L formalism, now in an
extended form with spectral parameter. For the proof of the second statement
we use the classical Lax pair. As a corollary of Theorem
\ref{Prop:spectral_invariance} we deduce that the mKdV flow in $\mathcal{O}%
_{\delta}(I\times{\mathbb{R}})$ with $\delta\le1/2$ preserves the spectrum of
the impedance operator (see Theorem \ref{Prop:spectral_invariance*}). In
addition, we prove Theorem \ref{Th:generalized_eigenfunctions} which states
that the evolution corresponding to the linear first-order differential
operator $Q_{\lambda}(s):=(4\lambda+2q(s,x))\partial_{x}-q_{x}(s,x)$, where
$s$ is a real parameter between $t$ and $t^{\prime}$ and $\lambda$ is the
spectral parameter, transforms any complete orthonormal system of generalized
eigenfunctions of the operator $L(t)$ to a complete orthonormal system of
generalized eigenfunctions of the operator $L(t^{\prime})$. Moreover, by
Corollary \ref{Prop:evolution}, the solution of the evolution equation
involving the third-order differential operator $A(t):=-4\partial_{x}%
^{3}+6q(t,x)\partial_{x}+3q_{x}(t,x)$ appearing in the classical Lax pair for
KdV, can be obtained in terms of the solution of the first-order evolution
equation $\psi_{t}=Q_{\lambda}(t)\psi$ with the spectral parameter $\lambda
\in\mathbb{R}$.

Growing solutions of evolution equations such as KdV and mKdV require the
development of new techniques for their study and are of interest by
themselves. They recently attracted a lot of attention. In \cite{Dub},
Dubrovin studied Hamiltonian perturbations of the (simplest) hyperbolic
equation $u_{t}+a(u)u_{x}=0$ in one space dimension. He conjectured that the
behavior of a solution to the perturbed equation near a point where the
gradient of the corresponding solution of the unperturbed equation blows up,
is \emph{universal}. This means that the behavior is (essentially) independent
of the choice of the (generic) Hamiltonian perturbation and of the (generic)
solution of the perturbed equation. In fact, he conjectured that the behavior
of solutions of the perturbed equations near such points is described by a
special smooth globally (in $X,T$) defined solution $U(X,T)$ of an integrable
fourth-order ODE in the variable $X$ which depends on a (real) parameter $T$.
When viewed as function of $X$ and $T$, $U(X,T)$ satisfies the KdV equation
and grows for $X\to\pm\infty$ as $\mp(6|X|)^{1/3}$.

\vspace{0.1cm}

\emph{Related work:} Beside the works \cite{B1,BS1,BS2,BS3}, we would like to
mention earlier work on unbounded solutions of KdV by Menikoff \cite{M} as
well as work of Kenig, Ponce, and Vega \cite{KPV}. Menikoff showed that for
initial data in $o_{1}(\mathbb{R})$, KdV can be solved in $C^{\infty
}(\mathbb{R}\times\mathbb{R})$ whereas Kenig, Ponce, and Vega studied
solutions of KdV in special classes of unbounded functions, different from the
ones considered in this paper. It was pointed out in \cite{M} that the KdV
flow with initial data in $o_{1}(\mathbb{R})$ preserves the discrete spectrum
of the Schr\"{o}dinger operator $L(t)$. We remark that the Miura map has been
used previously to obtain solutions of mKdV from solutions of KdV. In
particular, we mention the paper \cite{KT} where periodic solutions of low
regularity are obtained, and work of Gesztesy--Simon \cite{GS} and
Gesztesy--Schweiger--Simon \cite{GSS} for bounded solutions of mKdV. They use
techniques similar to those used here to construct solutions of the form
(\ref{representation}) to the mKdV\ equation from solutions $q(t,x)$ of the
KdV-equation with $q(t,x)$ and $q_{x}(t,x)$ bounded in $x$ for any
$t\in\mathbb{R}$. In contrast, we need to consider solutions of KdV\ which
grow at infinity and to derive precise asymptotics for the solutions
$\psi(t,x)$ of $-\psi_{xx}(t,x)+q(t,x)\psi(t,x)=0$. These asymptotics are
obtained by studying equation (\ref{evol1}).

Equations of the type (\ref{ldot}) have also appeared in the context of $2+1$
dimensionsal KdV-type equations. Manakov \cite{Man} observed that these latter
equations admit a representation of the form
\[
\frac{d}{dt}\left(  L_{2}-\lambda\right)  =\left[  L_{2}-\lambda,A\right]
+B(L_{2}-\lambda)
\]
where $L_{2}$, $A$, and $B$ are a so-called Manakov triple. See \cite{VN},
\cite{GN1}, \cite{GN2}, and further references in \cite{GN2} for details.

\emph{Acknowledgments:} We would like to thank A. Its and R. G. Novikov for
useful comments on an earlier version of our paper as well as B. Dubrovin who
introduced us to his recent work.

\section{Preliminaries}\label{Sec:auxiliaries} 
In this section, for the convenience of the reader, we describe some properties
of Q-L pairs which are mostly well-known.

Suppose that $q\in C^{\infty}(\mathbb{R}\times\mathbb{R})$ and consider
the differential operators $Q(t)$ and $L(t)$ given by (\ref{opQ}) and (\ref{opL}), respectively. 
\begin{Lemma}
\label{Lem:Commutation} The operators $Q$ and $L$ satisfy the following commutator relation
\begin{equation}
{\dot{L}}=[Q,L]+4q_{x}L+KdV(q)
\end{equation}
where ${\dot{L}}=q_{t}(t,x)$ and $KdV(q)=q_{t}-6qq_{x}+q_{xxx}$. In
particular,
\begin{equation}
KdV(q)=0\,\,\text{ iff }\,\,{\dot{L}}=[Q,L]+4q_{x}L.
\end{equation}
\end{Lemma}
The proof of the lemma is straightforward.

We remark that the operator $Q$ results from formally dividing the operator
$A=-4\partial_{x}^{3}+6q_{x}\partial_{x}+3q_{x}$ by $L=-\partial_{x}^{2}+q$.
In fact, $A=4\partial_{x}L+Q$. A similar division of $A$ by $L-\lambda$
results in the operator $Q_{\lambda}$ considered in section 6 below.

Assume that $q\in C^{\infty}(\mathbb{R}\times\mathbb{R})$ satisfies the KdV
equation and that for any $T>0$ there exists a constant $C_{T}>0$ such that
for any $|x|\ge1$ and $t\in[-T,T]$
\begin{equation}
\label{growth-condition2}|q(t,x)|\le C_{T}|x|\;.
\end{equation}
Let $\psi_{0}\in C^{\infty}(\mathbb{R})$ be an eigenfunction of $L(0)$ with
eigenvalue $0$, i.e.
\begin{equation}
\label{psi_0}L(0)\psi_{0}=0\;.
\end{equation}
Consider the equation
\begin{align}
\psi_{t}(t,x)  &  =Q(t)\psi(t,x)\label{evolution1}\\
\psi|_{t=0}  &  =\psi_{0} . \label{evolution2}%
\end{align}
By Lemma \ref{Lem:EDE}, the initial value problem (\ref{evolution1}%
)-(\ref{evolution2}) has a unique solution $\psi(t,x)$ in $C^{\infty
}(\mathbb{R}\times\mathbb{R})$.

\begin{Prop}
\hspace{-2mm}\textbf{.}\label{Prop:eigenfunction-transport} If $q\in
C^{\infty}(\mathbb{R}\times\mathbb{R})$ is a solution of KdV satisfying the
growth condition (\ref{growth-condition2}), and $\psi(t,x) \in C^{\infty
}(\mathbb{R}\times\mathbb{R})$ solves (\ref{evolution1})-(\ref{evolution2}),
then
\begin{equation}
\label{Schrodinger}L(t)\psi(t,x)=0\;\;\forall t,x\in\mathbb{R}.
\end{equation}
If, in addition, $\psi_{0}(x)>0$ $\forall x\in\mathbb{R},$ then $\psi(t,x)>0$
$\forall x,t\in\mathbb{R}$\;.
\end{Prop}

\emph{Proof.} Let $\varphi(t,x):=L(t)\psi(t,x)$. It follows from (\ref{psi_0})
that $\varphi|_{t=0}=0$. Using Lemma \ref{Lem:Commutation} and
(\ref{evolution1}) one obtains
\begin{align}
\varphi_{t}  &  ={\dot L}\psi+L\psi_{t}\nonumber\\
&  =([Q,L]+4q_{x}L)\psi+LQ\psi\nonumber\\
&  =Q(L\psi)+4q_{x}(L\psi)\nonumber\\
&  =2q\varphi_{x}+3q_{x}\varphi\;. \label{EPEonP}%
\end{align}
Hence $\varphi=\varphi(t,x)$ is a solution of the initial value problem
\begin{align}
\varphi_{t}(t,x)  &  =2q(t,x)\varphi_{x}(t,x)+3q_{x}(t,x)\varphi
(t,x)\nonumber\\
\varphi|_{t=0}  &  =0\;.\nonumber
\end{align}
Applying Lemma \ref{Lem:EDE} again, we obtain that $\varphi\equiv0$.

The last statement of the proposition follows immediately from claim $(b)$ of
Lemma \ref{Lem:EDE}. \hfill$\Box$ \vspace{3mm}

\begin{Prop}
\hspace{-2mm}\textbf{.}\label{Prop:eigenfunction_transport*} Assume that $r\in
C^{\infty}(\mathbb{R}\times\mathbb{R})$ is a solution of the initial value
problem (\ref{mKdV1})-(\ref{mKdV2}) for the mKdV equation and define
\begin{equation}
\rho(t,x):=\rho_{0}(t)e^{\int_{0}^{x}r(t,s)\;ds}%
\end{equation}
with normalizing factor $\rho_{0}(t)$ given by
\begin{equation}
\rho_{0}(t):=e^{\int_{0}^{t}(2r^{3}-r_{xx})|_{(\tau,0)}\;d\tau}\,.
\end{equation}
Then $\psi(t,x):=\rho(t,x)$ is a solution of (\ref{evolution1}%
)-(\ref{evolution2}), where $Q(t)=2q\partial_{x}-q_{x}$, $q=r_{x}+r^{2}$, and
$\psi_{0}(x):=e^{\int_{0}^{x}r_{0}(s)\;ds}$. If, in addition, $q=r_{x}+r^{2}$
satisfies the growth condition (\ref{growth-condition2}), then $\rho(t,x)$ is
the unique solution of (\ref{evolution1})-(\ref{evolution2}) in $C^{\infty
}(\mathbb{R}\times\mathbb{R})$.
\end{Prop}

\emph{Proof of Proposition \ref{Prop:eigenfunction_transport*}.} Using that
$q=r_{x}+r^{2},$ one easily sees that $\rho(t,x)$ satisfies the equation
$L(t)\rho=0$. Differentiating the latter identity with respect to $t$ and
using Lemma \ref{Lem:Commutation} together with the fact that $q=r_{x}+r^{2}$
satisfies KdV (cf. \cite{Miura}), we obtain
\begin{align}
0  &  ={\dot L}\rho+L\rho_{t}\nonumber\\
&  =([Q,L]+4q_{x}L)\rho+L\rho_{t}.\nonumber
\end{align}
Using the fact that $L(t)\rho=0$, one then gets
\begin{equation}
0 = - L(Q\rho)+L\rho_{t}=L(\rho_{t}-Q\rho).
\end{equation}
Hence, with $f(t,x):=\rho_{t}-Q\rho$ one has for any $t,x\in\mathbb{R}$
\begin{equation}
\label{secondorder}-f_{xx}(t,x)+q(t,x)f(t,x)=0.
\end{equation}
A direct computation shows that
\begin{equation}
\label{secondorder-data}f(t,0)=0\;\mbox{and}\;f_{x}(t,0)=0\;\forall
t\in\mathbb{R}.
\end{equation}
By the uniqueness of the solutions of (\ref{secondorder}%
)-(\ref{secondorder-data}) for any fixed $t\in\mathbb{R}$, we conclude that
$f(t,x)\equiv0$, and therefore $\rho_{t}=Q\rho$. The uniqueness of the
solution $\rho$ follows from Lemma \ref{Lem:EDE} together with the assumption
that $q(t,x)$ satisfies the growth condition (\ref{growth-condition2}).
\hfill$\Box$ \vspace{3mm}

\begin{Coro}
\hspace{-2mm}\textbf{.} Assume that $q\in C^{\infty}(\mathbb{R}\times
\mathbb{R})$ solves the KdV equation and satisfies the growth condition
(\ref{growth-condition2}). Let $\phi,\psi\in C^{\infty}(\mathbb{R}%
\times\mathbb{R})$ be two solutions of (\ref{evolution1}) with initial data
$\phi|_{t=0}=\phi_{0}$ and $\psi|_{t=0}=\psi_{0}$ respectively where
$L(0)\phi_{0}=0$ and $L(0)\psi_{0}=0$. Then the Wronskian $W(\phi,\psi
):=\phi\psi_{x}-\psi\phi_{x}$ is independent of $t,x\in\mathbb{R}$.
\end{Coro}

\emph{Proof.} As $\phi(t,x)$ and $\psi(t,x)$ satisfy (\ref{Schrodinger}) (see
Proposition \ref{Prop:eigenfunction-transport}) we get that the Wronskian $W$
is independent of $x\in\mathbb{R}$. Using that $\phi_{xx}=q\phi$ and
$\psi_{xx}=q\psi$ one obtains
\begin{align}
W_{t}  &  =\phi_{t}\psi_{x}+\phi(\psi_{t})_{x}-\psi_{t}\phi_{x}-\psi(\phi
_{t})_{x}\nonumber\\
&  =(2q\phi_{x}-q_{x}\phi)\psi_{x}+\phi(2q\psi_{x}-q_{x}\psi)_{x}-(2q\psi
_{x}-q_{x}\psi)\phi_{x}\nonumber\\
&  -\psi(2q\phi_{x}-q_{x}\phi)_{x}\nonumber\\
&  =0\;.\nonumber
\end{align}
\hfill$\Box$ \vspace{3mm}

\begin{Th}
\hspace{-2mm}\textbf{.}\label{Th:mKdV} Consider the initial value problem
(\ref{mKdV1})-(\ref{mKdV2}) for the mKdV equation with smooth initial data
$r_{0}\in C^{\infty}(\mathbb{R})$. Suppose that the solution $q=q(t,x)$ of the
KdV equation (\ref{KdV}) with the initial data $q|_{t=0}=q_{0}:=r_{0}^{\prime
}+r_{0}^{2}$ is defined globally in time, $q\in C^{\infty}(\mathbb{R\times
\mathbb{R}})$, and satisfies the growth condition (\ref{growth-condition2}). Then

\begin{itemize}
\item[(a)] the evolution equation (\ref{evolution1})-(\ref{evolution2}) has a
unique, globally defined, positive solution $\psi(t,x)>0$ and the function
$r(t,x)=\psi_{x}(t,x)/\psi(t,x)$ is a global solution of the mKdV initial
value problem (\ref{mKdV1})-(\ref{mKdV2});

\item[(b)] if $r_{1},r_{2}\in C^{\infty}(\mathbb{R}\times\mathbb{R})$ are
solutions of the initial value problem of mKdV (\ref{mKdV1})-(\ref{mKdV2})
both having $q$ as their image with respect to the Miura map $r\mapsto
r_{x}+r^{2}$ (i.e., $\forall t,x\in\mathbb{R}$, ${r_{1}}_{x}(t,x)+r_{1}%
^{2}(t,x)={r_{2}}_{x}(t,x)+r_{2}^{2}(t,x)$), then $r_{1}\equiv r_{2}$.
\end{itemize}
\end{Th}

\begin{Rem}
\hspace{-2mm}\textbf{.} Loosely speaking, statement $(b)$ of Theorem
\ref{Th:mKdV} says that whenever KdV has a unique solution within a certain
class then mKdV has a unique solution within the corresponding class defined
by the Miura map.
\end{Rem}

\emph{Proof of Theorem \ref{Th:mKdV}.} $(a)$ Introduce
\begin{equation}
\label{psi_0-def}\psi_{0}(x)=e^{\int_{0}^{x}r_{0}(s)\;ds}.
\end{equation}
Clearly, $\psi_{0}(x)>0$ $\forall x\in\mathbb{R}$. As $q_{0}=r_{0}^{\prime
}+r_{0}^{2}$ one obtains from (\ref{psi_0-def}) that $L(0)\psi_{0}=0$. By
Proposition \ref{Prop:eigenfunction-transport}, the solution $\psi(t,x)$ of
(\ref{evolution1})-(\ref{evolution2}) in $C^{\infty}(\mathbb{R}\times
\mathbb{R})$ satisfies $L(t)\psi(t,x)=0$ $\forall t,x\in\mathbb{R}$. Moreover,
$\psi(t,x)>0$ $\forall t,x\in\mathbb{R}$. Consider the smooth function
$r(t,x)$ given by (\ref{representation}). It follows from (\ref{psi_0-def})
that $r|_{t=0}=r_{0}$. Taking into account that $L(t)\psi(t,x)=0$ one proves
by a straightforward calculation that
\begin{align}
mKdV(r)  &  :=r_{t}-6r^{2}r_{x}+r_{xxx}\nonumber\\
&  =-(\psi_{t}\psi_{x}-\psi\psi_{xt}-6q\psi_{x}^{2}+ 3\psi_{xx}^{2}+4\psi
_{x}\psi_{xxx}-\psi\psi_{xxxx})/\psi^{2}\nonumber
\end{align}
(See also formula (7.42) in \cite{GSS}.) Using that $\psi_{t}=Q\psi$ one gets
that $mKdV(r)=0$. This proves claim $(a)$.

Claim $(b)$ follows from Proposition \ref{Prop:eigenfunction_transport*}, as
the two solutions $r_{1},r_{2}$ lead to the same operator $Q$ (cf.
(\ref{opQ})) and the same initial data $\psi_{0}$ (cf. (\ref{evolution2})).
Indeed, as $r_{1}$ and $r_{2}$ are solutions of (\ref{mKdV1})-(\ref{mKdV2})
and $q=r_{1x}+r_{1}^{2}=r_{2x}+r_{2}^{2}$ we get from Proposition
\ref{Prop:eigenfunction_transport*} that for $k=1,2$,
\[
\rho_{k}(t,x):=\rho_{k,0}(t)e^{\int_{0}^{x}r_{k}(t,s)\;ds}\;\;\;{\text{ with
}}\;\;\; \rho_{k,0}(t):=e^{\int_{0}^{t}(2r_{k}^{3}-(r_{k})_{xx})|_{(\tau
,0)}\;d\tau}
\]
are solutions of the linear initial value problem (\ref{evolution1}%
)-(\ref{evolution2}) with the same initial data $\psi_{0}(x)=\rho
_{k}(0,x)=e^{\int_{0}^{x}r_{0}(s)\;ds}$. As $q$ satisfies the growth condition
(\ref{growth-condition2}) the solution of (\ref{evolution1})-(\ref{evolution2}%
) is unique and therefore $\rho_{1}\equiv\rho_{2}$. In particular,
$r_{1}=\frac{\rho_{1x}}{\rho_{1}}=\frac{\rho_{2x}}{\rho_{2}}=r_{2}$.
\hfill$\Box$ \vspace{3mm}

\section{Proof of Theorem \ref{Th:mKdV-S}}

\label{Sect:mKdV-S} The purpose of this section is to prove Theorem
\ref{Th:mKdV-S}. In the sequel we will need the classes
\[
\mathcal{S}_{\beta+1}^{*}(\mathbb{R}\times\mathbb{R}):=\{f\in C^{\infty
}(\mathbb{R}\times\mathbb{R})\,|\, f_{x}\in\mathcal{S}_{\beta}(\mathbb{R}%
\times\mathbb{R})\}\,.
\]
where $\beta$ is a given real number. Note that the operator of integration,
$f(t,x) \mapsto\int_{0}^{x}f(t,s) ds,$ maps $\mathcal{S}_{\beta}%
(\mathbb{R}\times\mathbb{R})$ to $\mathcal{S}_{\beta+1}^{*}(\mathbb{R}%
\times\mathbb{R})$ whereas the operator of differentiation, $f(t,x)
\mapsto\partial_{x}f(t,x)$, maps $\mathcal{S}_{\beta+1}^{*}(\mathbb{R}%
\times\mathbb{R})$ to $\mathcal{S}_{\beta}(\mathbb{R}\times\mathbb{R})$ for
any $\beta\in\mathbb{R}$. Analogously one defines $\mathcal{S}_{\beta+1}%
^{*}(\mathbb{R})$.

The following Lemma describes the functions from $\mathcal{S}_{\beta+1}%
^{*}(\mathbb{R}\times\mathbb{R})$ in terms of their asymptotics at $\pm\infty$.

\begin{Lemma}
\hspace{-2mm}\textbf{.}\label{Prop:S*description} $f\in\mathcal{S}_{\beta
+1}^{*}(\mathbb{R}\times\mathbb{R})$ if and only if $f\in C^{\infty
}(\mathbb{R}\times\mathbb{R})$ and it has an asymptotic expansion for $x\to
\pm\infty$ of the form
\begin{equation}
\label{e:S*asymptotics}f(t,x)\sim\left\{
\begin{array}
[c]{lcl}%
\sum_{k=0}^{\infty}\;a_{k}^{\pm}(t)\,(\pm x)^{\beta_{k}+1}+a_{*}^{\pm}%
(t)\log(\pm x) & \mbox{\rm if} & \beta+1\ge0\\
&  & \\
c^{\pm}(t)+\sum_{k=0}^{\infty}\; a_{k}^{\pm}(t)\,(\pm x)^{\beta_{k}+1} &
\mbox{\rm if} & \beta+1<0
\end{array}
\right.
\end{equation}
where $\beta=\beta_{0}>\beta_{1}>...$ with $\lim\limits_{k\to\infty}\beta
_{k}=-\infty$ and $a_{k}^{\pm}$, $a_{*}^{\pm}$, and $c_{\pm}$ are functions of
$t$ in $C^{\infty}(\mathbb{R})$. The same result holds in $\mathcal{S}%
_{\beta+1}^{*}(\mathbb{R})$.
\end{Lemma}

In particular, if $\beta+1\ge0$ then the leading term of the asymptotic
expansion of $f$ is $a_{0}^{\pm}(t)(\pm x)^{\beta+1}$ (for $\beta>-1$) or
$a_{*}^{\pm}(t)\log(\pm x)$ (for $\beta=-1$). If $\beta+1<0$, the leading term
is $c^{\pm}(t)$ followed by $a_{0}^{\pm}(t)(\pm x)^{\beta+1}$. The asymptotic
relations should be understood similarly to \eqref{at+infty},
\eqref{at-infty}. For example, the first relation in \eqref{e:S*asymptotics}
means that for any compact interval $J\subseteq\mathbb{R}$, $i,j\ge0$, and any
$N\ge0$ with $\beta_{N}+1<0$, there exists a constant $C_{J,N,i,j}>0$ such
that for any $|x|\ge1$ and any $t\in J$
\begin{equation}
\label{e:asymptotics*}\Big|\;\partial_{t}^{i}\partial_{x}^{j}%
\Big(f(t,x)-\Big(a_{*}^{\pm}(t)\log(\pm x)+ \sum_{k=0}^{N} a_{k}^{\pm}(t)(\pm
x)^{\beta_{k}+1}\Big)\Big)\Big|\le C_{J,N,i,j}|x|^{(\beta_{N+1}+1)-j}.
\end{equation}
\emph{Proof of Lemma \ref{Prop:S*description}.} If $f\in C^{\infty}%
(\mathbb{R}\times\mathbb{R})$ has an asymptotic expansion as in
\eqref{e:S*asymptotics}, then clearly, $f_{x}\in\mathcal{S}_{\beta}%
(\mathbb{R}\times\mathbb{R})$, hence $f\in\mathcal{S}_{\beta+1}^{*}%
(\mathbb{R}\times\mathbb{R})$. Let us prove the converse statement. As the
asymptotic expansions for $x\to+\infty$ and $x\to-\infty$ of an element
$f\in\mathcal{S}_{\beta+1}^{*}(\mathbb{R}\times\mathbb{R})$ are obtained in a
similar way let us consider the case $x\to+\infty$ only. First we treat the
case where $\beta+1\ge0$. By definition, for an element $f\in\mathcal{S}%
_{\beta+1}^{*}(\mathbb{R}\times\mathbb{R})$, $f_{x}\in\mathcal{S}_{\beta
}(\mathbb{R}\times\mathbb{R})$ and hence has an asymptotic expansion
\begin{equation}
\label{f_xasymptotics}f_{x}\sim\sum_{k=0}^{\infty}b_{k}^{+}(t)\,x^{\beta_{k}%
}\;\;\;\mbox{as}\;\;\;x\to\infty
\end{equation}
where $\beta=\beta_{0}>\beta_{1}>...$ with $\lim\limits_{k\to\infty}\beta
_{k}=\infty$. Without loss of generality we assume that $\beta_{m}=-1$ for
some $m\ge0$.\footnote{Take $b_{m}^{+}(t)\equiv0$ if necessary.} Formally, the
claimed result is obtained by integrating term by term the right hand side of
\eqref{f_xasymptotics} with respect to the $x$-variable. In order to make this
argument rigorous we argue as follows: For any $N\ge m+1$ and $x\in\mathbb{R}$
consider the quantity
\begin{equation}
\label{e:Q_N}Q_{N}(t,x):=-\chi_{+}(x)\int_{x}^{\infty}\Big(f_{x}%
(t,s)-\sum_{k=0}^{N} b_{k}^{+}(t) s^{\beta_{k}}\Big)\,ds
\end{equation}
where $\chi_{+}(x)$ is a smooth cut-off function with $\chi_{+}(x)=0$ for
$x\le1/2$ and $\chi_{+}(x)=1$ for $x\ge1$. As $f_{x}\in\mathcal{S}_{\beta
}(\mathbb{R}\times\mathbb{R})$ and $\beta_{m+1}<-1$ it follows that the
improper integral in \eqref{e:Q_N} exists and if $x\ge1$, $\partial_{x}%
Q_{N}(t,x)=f_{x}(t,x)-\sum_{k=0}^{N} b_{k}^{+}(t) x^{\beta_{k}}$. Hence,
$\partial_{x} Q_{N}$ is in $\mathcal{S}_{\beta_{N+1}}(\mathbb{R}%
\times\mathbb{R})$. We claim that $Q_{N}$ is in $\mathcal{S}_{\beta_{N+1}%
+1}(\mathbb{R}\times\mathbb{R})$. To show this, it remains to estimate
$\partial_{t}^{i}Q_{N}(t,x)$. It follows from \eqref{f_xasymptotics} that for
any compact interval $J\subseteq\mathbb{R}$, $i\ge0$, and $N\ge m+1$, there
exists a constant $C_{J,N,i}>0$ such that for any $x\ge1$, $t\in J$
\begin{align}
|\partial_{t}^{i} Q_{N}(t,x)|  &  \le\int_{x}^{\infty}|\partial_{t}^{i}%
(f_{x}(t,s)-\sum_{k=0}^{N} b_{k}^{+}(t) s^{\beta_{k}})|\,ds\nonumber\\
&  \le C_{J,N,i}\int_{x}^{\infty}s^{\beta_{N+1}}\,ds\nonumber\\
&  \le C_{J,N,i}\,\frac{x^{\beta_{N+1}+1}}{|\beta_{N+1}+1|} \,. \label{estQ_N}%
\end{align}
Computing the integral in \eqref{e:Q_N} one gets for $x\ge1$
\begin{align}
\label{Q_N}Q_{N}(t,x)=f(t,x)-\Big(c^{+}(t)+b_{m}^{+}(t)\log x+\sum_{0\le k\le
N,k\ne m}\frac{b_{k}^{+}(t)}{\beta_{k}+1}\,x^{\beta_{k}+1}\Big)
\end{align}
where
\begin{equation}
\label{c^+}c^{+}(t):=f(t,1)-\sum_{k=0}^{m-1}\frac{b_{k}^{+}(t)}{\beta_{k}+1}
+\int_{1}^{\infty}\Big(f_{x}(t,s)-\sum\limits_{k=0}^{m} b_{k}^{+}(t)
s^{\beta_{k}}\Big)\,ds\,.
\end{equation}
(Note that the integral in \eqref{c^+} converges as the integrand is estimated
locally uniformly in $t$ by $O(s^{\beta_{m+1}})$ with $\beta_{m+1}<-1$.) The
desired estimate \eqref{e:asymptotics*} of $f(t,x)$ for $x\to+\infty$ follows
from \eqref{f_xasymptotics}, \eqref{estQ_N}, and \eqref{Q_N}.

The case $\beta+1<0$ is treated in a similar way. Actually, it is easier than
the case $\beta+1\ge0$. \hfill$\Box$ \vspace{3mm}

\noindent\emph{Proof of Theorem \ref{Th:mKdV-S}.} We will show that the
claimed results follow from Theorem \ref{Th:mKdV} and Lemma
\ref{Lem:EDE-asymptotics} stated below, combined with results in \cite{BS1,B1}
- see Appendix B for a summary of these results. Indeed, for a given $r_{0}%
\in\mathcal{S}_{\beta}(\mathbb{R}),$ the Miura image $q_{0}:={r_{0}}_{x}%
+r_{0}^{2}$ belongs to $\mathcal{S}_{\delta}(\mathbb{R})$ with $\delta
:=\max\{2\beta,\beta-1\}<1$. According to the results in \cite{BS1,B1} (cf.
Theorem \ref{Th:SB-as}, \ref{Th:SB-uniqueness} in Appendix B) there exists a
unique solution $q\in\mathcal{S}_{\delta}(\mathbb{R}\times\mathbb{R})$ of the
KdV equation (\ref{KdV}) with initial data $q|_{t=0}=q_{0}$. As $\delta<1$ the
solution $q=q(t,x)$ satisfies the growth condition (\ref{growth-condition2}).
In particular, according to Proposition \ref{Prop:eigenfunction-transport} the
linear initial value problem
\begin{align}
\psi_{t}(t,x)  &  =2q(t,x)\psi_{x}(t,x)-q_{x}(t,x)\psi(t,x)\label{psiEDE1}\\
\psi|_{t=0}  &  =\psi_{0}(x):=e^{\int_{0}^{x}r_{0}(s)\;ds} \label{psiEDE2}%
\end{align}
has a unique $C^{\infty}$-solution. This solution $\psi$ is defined globally
in time. It is strictly positive everywhere, and satisfies $-\psi_{xx}%
+q\psi=0$ $\forall t,x\in\mathbb{R}$. It follows from item $(a)$ of Theorem
\ref{Th:mKdV} that the function $r(t,x)=\psi_{x}(t,x)/\psi(t,x)$ is a solution
of (\ref{mKdV1})-(\ref{mKdV2}). It is easy to see that the function
\[
p=p(t,x):=\log\psi(t,x)
\]
satisfies
\begin{align}
p_{t}(t,x)  &  =2q(t,x)p_{x}(t,x)-q_{x}(t,x)\label{EDE-asymptotics1}\\
p|_{t=0}  &  =p_{0}(x) \label{EDE-asymptotics2}%
\end{align}
with initial data $p_{0}(x)=\int_{0}^{x}r_{0}(s)\;ds\in\mathcal{S}_{\beta
+1}^{*}(\mathbb{R})$. According to Lemma \ref{Lem:EDE-asymptotics} below the
function $p(t,x)$ belongs to $\mathcal{S}_{\beta+1}^{*}(\mathbb{R}%
\times\mathbb{R})$ and therefore $r=\partial_{x} p\in\mathcal{S}_{\beta
}(\mathbb{R}\times\mathbb{R})$.

The uniqueness of the solution $r=r(t,x)$ in the class $\mathcal{S}_{\beta
}(\mathbb{R}\times\mathbb{R})$ follows from the uniqueness of the solution
$q=q(t,x)$ in the class $\mathcal{S}_{\delta}(\mathbb{R}\times\mathbb{R})$
(cf. Theorem \ref{Th:SB-uniqueness} in Appendix B) and Theorem \ref{Th:mKdV}
$(b)$.\hfill$\Box$ \vspace{3mm}

The proof of Theorem \ref{Th:mKdV-S} used the following lemma.

\begin{Lemma}
\hspace{-2mm}\textbf{.}\label{Lem:EDE-asymptotics} Let $\beta<1/2$ and
$\delta:=\max\{2\beta,\beta-1\}<1$. Consider the initial value problem
(\ref{EDE-asymptotics1})-(\ref{EDE-asymptotics2}) where $q\in\mathcal{S}%
_{\delta}(\mathbb{R}\times\mathbb{R})$. Then for any initial data $p_{0}%
\in\mathcal{S}_{\beta+1}^{*}(\mathbb{R})$ there exists a solution
$p\in\mathcal{S}_{\beta+1}^{*}(\mathbb{R}\times\mathbb{R})$ of
(\ref{EDE-asymptotics1})-(\ref{EDE-asymptotics2}). This solution is unique in
$C^{\infty}(\mathbb{R}\times\mathbb{R})$.
\end{Lemma}

In order to prove Lemma \ref{Lem:EDE-asymptotics} we will first construct
formal series $\chi_{\pm}(t,x)$ having the form \eqref{e:S*asymptotics} and
satisfying the evolution equation
\eqref{EDE-asymptotics1}-\eqref{EDE-asymptotics2} formally for $x\to\pm\infty
$. As the cases $x\to+\infty$ and $x\to-\infty$ are treated in the same way we
restrict our attention only to the case $x\to+\infty$.

Let $p_{0} \in\mathcal{S}_{\beta+1}^{*}(\mathbb{R})$ be the initial data in
\eqref{EDE-asymptotics1}-\eqref{EDE-asymptotics2}. By Lemma
\ref{Prop:S*description}, $p_{0}(x)$ has an asymptotic expansion for
$x\to+\infty$ of the form
\[
p_{0}(x)\sim\left\{
\begin{array}
[c]{lcl}%
\sum_{k=0}^{\infty}p_{k}^{+} x^{\beta_{k}+1}+p_{*}^{+}\log{x} &
\mbox{\rm if} & \beta+1\ge0\\
&  & \\
c^{+}+\sum_{k=0}^{\infty}p_{k}^{+} x^{\beta_{k}+1} & \mbox{\rm if} & \beta+1<0
\end{array}
\right.
\]
where $\beta_{0}:=\beta<1/2$ and $\beta_{0}>\beta_{1}>...$, $\lim
\limits_{k\to\infty}\beta_{k}\,=\,-\infty$. As a solution $p$ of
\eqref{EDE-asymptotics1}-\eqref{EDE-asymptotics2} gives rise to the
$1$-parameter family of solutions $p+\mathop{\tt const}$, we can assume
without loss of generality that the constant $c^{+}$ in the asymptotic
expansion for $p_{0}(x)$ vanishes,
\begin{equation}
\label{p_0-asymptotics}p_{0}(x)\sim\sum_{k=0}^{\infty}p_{k}^{+} x^{\beta
_{k}+1}+p_{*}^{+}\log{x}\;\;\mbox{as}\;\;x\to\infty\,.
\end{equation}
Here $\beta<1/2$ but it need not be true that $\beta+1\ge0$. By assumption,
$q\in\mathcal{S}_{\delta}(\mathbb{R}\times\mathbb{R})$ and hence it has an
asymptotic expansion for $x\to+\infty$ of the form
\begin{equation}
\label{q-asymptotics}q(t,x)\sim\sum_{k=0}^{\infty}c_{k}^{+}(t) x^{\delta_{k}}%
\end{equation}
where $\delta_{0}=\delta:=\max\{2\beta,\beta-1\}<1$ and $\delta_{0}>\delta
_{1}>...$, $\lim\limits_{k\to\infty}\delta_{k}=-\infty$. Consider the set
\[
\Delta:=\{\delta_{k}\}_{k\ge0}.
\]
In order to find a \emph{formal} solution $\chi_{+}(t,x)$ of
\eqref{EDE-asymptotics1}-\eqref{EDE-asymptotics2} we will have to extend the
set of exponents $\{\beta_{k}+1\}_{k\ge0}$ appearing in
\eqref{p_0-asymptotics} to a larger discrete set ${\bar B}$ with the same
upper limit as $\{\beta_{k}+1\}_{k\ge0}$ so that the exponents appearing in
the asymptotic expansions of the left- and right-hand side of
\eqref{EDE-asymptotics1} belong to ${\bar B}$. To construct ${\bar B}$ we
first need to extend the set $\Delta$.

\begin{Lemma}
\hspace{-2mm}\textbf{.}\label{Lem:B} There exists an unbounded discrete set
${\bar\Delta}\subseteq\mathbb{R}$ with $\Delta\subseteq{\bar\Delta}$ such that

\begin{itemize}
\item[(a)] $\max{\bar\Delta}=\max\Delta=\delta<1$;

\item[(b)] if $\delta^{\prime},\delta^{\prime\prime}\in{\bar\Delta}$ then
$\delta^{\prime}+\delta^{\prime\prime}-1\in{\bar\Delta}$;

\item[(c)] if $\delta^{\prime}\in{\bar\Delta}$ then $\delta^{\prime}-1\in
{\bar\Delta}$.
\end{itemize}
\end{Lemma}

\emph{Proof.} First note that a set ${\bar\Delta}$ satisfies (b) iff
${\bar\Delta}-1:=\{\delta^{\prime}-1\,|\,\delta^{\prime}\in{\bar\Delta}\}$
satisfies
\begin{equation}
\label{b}\delta^{\prime},\delta^{\prime\prime}\in{\bar\Delta}%
-1\;\;\;\mbox{implies}\;\;\;\delta^{\prime}+\delta^{\prime\prime}\in
{\bar\Delta}-1\;.
\end{equation}
It is easy to see that the set $\Delta_{1}\subseteq\mathbb{R}$,
\[
\Delta_{1}:=\Big\{\sum\limits_{i\in J}\delta_{i}\,|\,\delta_{i}\in\Delta-1,
J\subseteq\mathbb{Z}_{\ge0}\; \mbox{is finite and}\;J\ne\emptyset\Big\}
\]
is discrete, satisfies (\ref{b}) and that $\max\Delta_{1}=\delta-1$. Consider
the set
\[
{\bar\Delta}_{1}:=\{\delta^{\prime}-k\,|\,\delta^{\prime}\in\Delta_{1},
k\in\mathbb{Z}_{\ge0}\}.
\]
Then
\[
\delta^{\prime},\delta^{\prime\prime}\in{\bar\Delta}_{1}%
\;\;\;\mbox{implies}\;\;\;\delta^{\prime}+\delta^{\prime\prime}\in{\bar\Delta
}_{1};
\]
in addition, ${\bar\Delta}_{1}$ is unbounded and discrete. Moreover $\max
{\bar\Delta}_{1}=\delta-1$, and $\delta^{\prime}-1\in{\bar\Delta}_{1}$ for any
$\delta^{\prime}\in{\bar\Delta_{1}}$. Hence, the set ${\bar\Delta}%
:={\bar\Delta}_{1}+1$ satisfies claims (a)-(c) of the lemma. \hfill$\Box$
\vspace{3mm}

We extend in the sum in (\ref{q-asymptotics}) the set of exponents $\Delta$ to
${\bar\Delta}$ by setting the new coefficients in (\ref{q-asymptotics}) all
equal to zero. Hence without loss of generality, one can -- and in the sequel
we will -- assume that the set of the exponents $\Delta=\{\delta_{k}\}_{k\ge
0}$ in (\ref{q-asymptotics}) satisfies conditions (a)-(c) of Lemma \ref{Lem:B}.

Let us now introduce the following subsets of $\mathbb{R}$,
\[
B:=\{\beta_{k}\}_{k\ge0}
\]
and
\begin{equation}
{\bar B}:=\{\beta^{\prime}+\delta^{\prime}\,|\,\beta^{\prime}\in
B,\;\delta^{\prime}\in\Delta\}\cup\Delta\cup\{\beta^{\prime}+1\,|\,\beta
^{\prime}\in B\}
\end{equation}

\begin{Lemma}
\hspace{-2mm}\textbf{.}\label{Lem:barA} The set ${\bar B}$ is discrete and has
the following properties:

\begin{itemize}
\item[$(i)$] $\max{\bar B}=\beta+1$;

\item[$(ii)$] if $\delta^{\prime}\in\Delta$ and $\beta^{\prime}\in{\bar B}$,
then $\delta^{\prime}+\beta^{\prime}-1\in{\bar B}$;

\item[$(iii)$] the set $\{\delta^{\prime}-1\;|\:\delta^{\prime}\in\Delta\}$ is
contained in ${\bar B}$.
\end{itemize}
\end{Lemma}

\emph{Proof of Lemma \ref{Lem:barA}.} The proof that $\bar B$ is discrete
follows from the arguments used in the proof of Lemma \ref{Lem:B}.

\noindent$(i)$ As $\beta<1/2$ and $\delta=\max\{2\beta,\beta-1\}<1$ one gets
$\max{\bar B}=\max\{\delta,\beta+\delta,\beta+1\}=\beta+1$.

\noindent$(ii)$ follows from the fact that $\Delta$ has property $(b)$ of
Lemma \ref{Lem:B}. Indeed, as any $\beta^{\prime}\in{\bar B}$ can be written
in the form $\beta^{\prime}=\beta^{\prime\prime}+\delta^{\prime\prime}$
($\beta^{\prime\prime}\in B$, $\delta^{\prime\prime}\in\Delta$),
$\beta^{\prime}=\delta^{\prime\prime}$, or $\beta^{\prime}=\beta^{\prime
\prime}+1$ and as by Lemma \ref{Lem:B}(b), for any $\delta^{\prime}\in\Delta$,
one has $\delta^{\prime\prime\prime}:=\delta^{\prime}+\delta^{\prime\prime
}-1\in\Delta$, it follows that
\[
\delta^{\prime}+\beta^{\prime}-1=\left\{
\begin{array}
[c]{l}%
\delta^{\prime}+(\beta^{\prime\prime}+\delta^{\prime\prime})-1=\delta
^{\prime\prime\prime}+\beta^{\prime\prime}\in{\bar B}\\
\delta^{\prime}+\delta^{\prime\prime}-1\in\Delta\subseteq{\bar B}\\
\delta^{\prime}+(\beta^{\prime\prime}+1)-1=\delta^{\prime}+\beta^{\prime
\prime}\in{\bar B}%
\end{array}
\right.
\]
$(iii)$ It follows from statement (c) of Lemma \ref{Lem:B} that for any
$\delta^{\prime}\in\Delta,$ one has $\delta^{\prime}-1\in\Delta$ and as
$\Delta\subseteq{\bar B}$, $(c)$ then follows. \hfill$\Box$ \vspace{3mm}

\noindent\emph{Proof of Lemma \ref{Lem:EDE-asymptotics}.} First we prove that
for any
\begin{equation}
\label{q-asymptotic}q(t,x)\sim\sum_{k=0}^{\infty}c_{k}^{\pm}(t)(\pm
x)^{\delta_{k}}\;\;\mbox{as}\;\;x\to\pm\infty
\end{equation}
with exponents $\Delta=\{\delta_{k}\}_{k\ge0}$, $\delta_{0}=\delta>\delta
_{1}>...$, satisfying claims (a)-(b) of Lemma \ref{Lem:B}, the initial value
problem (\ref{EDE-asymptotics1})-(\ref{EDE-asymptotics2}) with $p_{0}(x)$
satisfying \eqref{p_0-asymptotics} has a formal solution $\chi_{+}(t,x)$ given
by ($t\in\mathbb{R}, x>0$)
\begin{equation}
\label{p-asymptotic}\chi_{+}(t,x)=\sum_{k=0}^{\infty}a_{k}^{+}(t)
x^{{\bar\beta}_{k}+1}+a_{*}^{+}(t)\log x
\end{equation}
where $\{{\bar\beta}_{k}+1\}_{k\ge0}={\bar B}$ with ${\bar\beta_{0}>\bar
\beta_{1}>...}$. The existence of a formal solution $\chi_{-}(t,x)$ for
$t\in\mathbb{R}$, $x<0$ follows by the same arguments. Let us stress that the
exponents $\beta_{k}+1$ in the asymptotic expansion of the initial data
$p_{0}(x)$ (cf. (\ref{p_0-asymptotics})) belong to the set $\{\beta^{\prime
}+1\,|\,\beta^{\prime}\in B\}$ which is included in the larger set $\bar B$.
Hence, the coefficients $a_{k}^{+}(0)$ in the asymptotic expansion
\eqref{p-asymptotic} evaluated at $t=0$ are zero or coincide with some of the
constants $p_{k}^{+}$. Moreover $a_{*}^{+}(0)=p_{*}^{+}$.

Substituting (\ref{q-asymptotic}) and (\ref{p-asymptotic}) into
(\ref{EDE-asymptotics1}) and using the notation $\empty^{\mathbf{\cdot}}%
=\frac{d}{dt}$, one obtains, in the case $x>0$,
\begin{align}
{\dot a}_{*}^{+}(t)\log x+ {\dot a}_{0}^{+}(t)x^{\beta+1}+{\dot a}_{1}%
^{+}(t)x^{{\bar\beta}_{1}+1}+{\dot a}_{2}^{+}(t)x^{{\bar\beta}_{2}+1}+
...=\nonumber\\
=2\Big(c_{0}^{+}(t)x^{\delta}+c_{1}^{+}(t)x^{\delta_{1}}+...\Big) \Big(a_{*}%
^{+}(t)x^{-1}+\nonumber\\
(\beta+1)a_{0}^{+}(t)x^{\beta}+({\bar\beta}_{1}+1)a_{1}^{+}(t)x^{{\bar\beta
}_{1}}+...\Big)\label{compareEq}\\
-\Big(\delta c_{0}^{+}(t)x^{\delta-1}+\delta_{1} c_{1}^{+}(t)x^{\delta_{1}%
-1}+...\Big)\nonumber
\end{align}

The maximal power of $x$ on the right side of (\ref{compareEq}) is not bigger
than $m_{r}=\max\{\beta+\delta,\delta-1\}$. As $\beta<1/2$ and $\delta
=\max\{2\beta,\beta-1\}<1$, one obtains that $\beta+1>m_{r}$. Hence, ${\dot
a}_{0}^{+}(t)=0$ and thus $a_{0}^{+}(t)=a_{0}^{+}(0)$. Comparing the
coefficients in (\ref{compareEq}) we also obtain that ${\dot a}_{*}^{+}(t)=0$
and hence $a_{*}^{+}(t)=a_{*}^{+}(0)$.

Comparing the coefficients in (\ref{compareEq}) corresponding to terms of
order ${\bar\beta}_{k}+1$ in $x$ one obtains that for any $k\ge1$
\begin{equation}
\label{theuptriangle}{\dot a}_{k}^{+}(t)=P_{k}^{+}(a_{0}^{+},a_{1}%
^{+}(t),...,a_{k-1}^{+}(t))+F_{k}^{+}(t)
\end{equation}
where $P_{k}^{+}$ is a linear combination of the variables $a_{0}%
^{+},...,a_{k-1}^{+}$ with coefficients which are smooth functions of
$t\in\mathbb{R}$. The term $F_{k}^{+}(t)$ is equal to $2c_{i_{k}}^{+}%
(t)p_{*}^{+} -\delta_{i_{k}}c_{i_{k}}^{+}(t)$ iff there exists an index
$i_{k}\ge0$ such that ${\bar\beta}_{k}+1=\delta_{i_{k}}-1$. If there is no
such $i_{k}$ then $F_{k}^{+}(t)\equiv0$. Let us prove formula
(\ref{theuptriangle}). It is clear from (\ref{compareEq}) that the right side
of (\ref{theuptriangle}) is a sum of a linear polynomial of the variables
$a_{0}^{+},a_{1}^{+},...$ and an inhomogeneous term $F_{k}^{+}(t)$ of the form
described above. Assume that there exists $a_{n}^{+}$, $n\ge k$, that enters
as a linear term on the right side of (\ref{theuptriangle}). Then clearly
there exists $m_{n}\ge0$ such that
\[
{\bar\beta}_{n}+\delta_{m_{n}}={\bar\beta}_{k}+1\;.
\]
As ${\bar\beta}_{n}\le{\bar\beta}_{k}$ and $\delta_{m_{n}}\le\delta<1$, it
follows that ${\bar\beta}_{n}+\delta_{m_{n}}<{\bar\beta}_{k}+1$. This
contradiction proves (\ref{theuptriangle}).

Integrating equation (\ref{theuptriangle}), we find the coefficients
$a_{k}^{+}(t)$ recursively in terms of the initial values $(a_{i}%
^{+}(0))_{0\le i\le k}$. Clearly, the formal solution $\chi_{+}(t,x)$
satisfies (\ref{compareEq}) and, by construction, $\chi_{+}(0,x)=p_{0}(x)$.
Arguing similarly we find a formal solution $\chi_{-}(t,x)$ for $t\in
\mathbb{R}$, $x<0$.

Next we show how the constructed formal solutions
\[
\chi_{\pm}(t,x)=\sum_{k=0}^{\infty}a_{k}^{\pm}(t)(\pm x)^{{\bar\beta}_{k}%
+1}+a_{*}^{\pm}(t)\log(\pm x)
\]
lead to a solution of (\ref{EDE-asymptotics1})-(\ref{EDE-asymptotics2}).
Choose $f(t,x)\in C^{\infty}(\mathbb{R}\times\mathbb{R})$ so that $f$ has
asymptotic expansions of the form
\begin{equation}
\label{formal_solution}f(t,x)\sim\sum_{k=0}^{\infty}a_{k}^{\pm}(t)(\pm
x)^{{\bar\beta}_{k}+1}+ a_{*}^{\pm}(t)\log(\pm x)\;\;\mbox{as}\;\;x\to
\pm\infty
\end{equation}
with coefficients $(a_{k}^{\pm}(t))_{k\ge0}$, $a^{\pm}_{*}(t)$ defined as
above. The existence of such a function $f$ follows, for example, from
\cite[Proposition 3.5]{Shubin}. Following \cite{BS1,B1} we will call the
function $f(t,x)$ an \emph{asymptotic solution} of (\ref{EDE-asymptotics1}%
)-(\ref{EDE-asymptotics2}). Let $f_{0}:=f|_{t=0}$.

With the help of the asymptotic solution $f$ we want to find a solution
$p(t,x)$ of (\ref{EDE-asymptotics1})-(\ref{EDE-asymptotics2}) of the form
\begin{equation}
\label{substitution}p(t,x):=f(t,x)+s(t,x)
\end{equation}
where $s(t,x)\in C^{\infty}(\mathbb{R}\times\mathbb{R})$ has to be determined
so that $s$ and all derivatives $\partial_{t}^{l}\partial_{x}^{k}s$ are fast
decaying as $|x|\to\infty$. Substituting (\ref{substitution}) into
(\ref{EDE-asymptotics1})-(\ref{EDE-asymptotics2}) one obtains the linear
evolution equation
\begin{align}
s_{t}(t,x)  &  =2q(t,x)s_{x}(t,x)+\eta(t,x)\label{EDE-S1}\\
s|_{t=0}  &  =s_{0}(x) \label{EDE-S2}%
\end{align}
where $\eta(t,x):=-f_{t}(t,x)+2q(t,x)f_{x}(t,x)-q_{x}(t,x)$ belongs to
$\mathcal{S}_{-\infty}(\mathbb{R}\times\mathbb{R})$ (as $f(t,x)$ is an
asymptotic solution of (\ref{EDE-asymptotics1})-(\ref{EDE-asymptotics2})) and
$s|_{t=0}=p_{0}(x)-f_{0}(x)\in\mathcal{S}(\mathbb{R})$, where as usual,
$\mathcal{S}(\mathbb{R})$ denotes the functions of Schwartz class. By
definition, $g\in C^{\infty}(\mathbb{R}\times\mathbb{R})$ belongs to the space
$\mathcal{S}_{-\infty}(\mathbb{R}\times\mathbb{R})$ iff for any compact
interval $J\subseteq\mathbb{R}$ and any $k,i,j\ge0$ there exists a constant
$C_{J,k,i,j}>0$ such that for any $|x|\ge1$ and $t\in J$
\[
|\partial^{i}_{t}\partial^{j}_{x}g(t,x)|\le C_{J,k,i,j}|x|^{-k}.
\]
In particular, if $g \in\mathcal{S}_{-\infty}(\mathbb{R}\times\mathbb{R})$
then for any given $t\in\mathbb{R},$ the function $g(t,\cdot)$ belongs to
$\mathcal{S}(\mathbb{R})$.

Due to Lemma \ref{Lem:EDE-S}, we can find a solution $s\in\mathcal{S}%
_{-\infty}(\mathbb{R}\times\mathbb{R})$ of \eqref{EDE-S1}-\eqref{EDE-S2} which
proves the existence part of Lemma \ref{Lem:EDE-asymptotics}. The uniqueness
of the solution $p(t,x)$ in $C^{\infty}(\mathbb{R}\times\mathbb{R})$ follows
from Lemma \ref{Lem:EDE}. \hfill$\Box$ \vspace{3mm}

\section{Proof of Theorem \ref{Th:mKdV-O} and Theorem \ref{Th:mKdV-o}}

\label{Sect:mKdV-Oo} In this section we prove the existence and the uniqueness
of solutions of the mKdV equation, globally in time, as stated in Theorem
\ref{Th:mKdV-O} and Theorem \ref{Th:mKdV-o}.

Before proving these theorems we introduce the following auxiliary spaces
\[
\mathcal{O}^{*}_{\beta+1}(\mathbb{R}\times\mathbb{R}):= \{f\in C^{\infty
}(\mathbb{R}\times\mathbb{R})\,|\,f_{x}\in\mathcal{O}_{\beta}(\mathbb{R}%
\times\mathbb{R})\}
\]
and
\[
{o}^{*}_{\beta+1}(\mathbb{R}\times\mathbb{R}):= \{f\in C^{\infty}%
(\mathbb{R}\times\mathbb{R})\,|\,f_{x}\in{o}_{\beta}(\mathbb{R}\times
\mathbb{R})\}\,.
\]
\noindent\emph{Proof of Theorem \ref{Th:mKdV-O}.} We follow the arguments in
the proof of Theorem \ref{Th:mKdV-S}. Given $r_{0}\in\mathcal{O}_{\beta
}(\mathbb{R})$ we get $q_{0}:=r_{0}^{\prime}+r^{2}_{0}$ which belongs to the
space $\mathcal{O}_{\delta}(\mathbb{R})$ with $\delta:=\max\{2\beta
,\beta-1\}<1$. By Theorem 2 in \cite{BS3} (cf. Theorem \ref{Th:SB-O}, Appendix
B) there exists a solution $q\in\mathcal{O}_{\delta}(\mathbb{R}\times
\mathbb{R})$ of the KdV initial value problem
\[
q_{t}-6qq_{t}+q_{xxx}=0,\;\;\;q|_{t=0}=q_{0}\;.
\]
As $\delta<1$ the solution $q(t,x)$ satisfies the growth condition
(\ref{growth-condition2}). Let $\psi(t,x)>0$ be the globally defined unique
solution of \eqref{psiEDE1}-\eqref{psiEDE2} (see also Proposition
\ref{Prop:eigenfunction-transport}). According to Theorem \ref{Th:mKdV} $(a)$
the function $r(t,x)=\psi_{x}(t,x)/\psi(t,x)$ is a solution of the mKdV
initial value problem (\ref{mKdV1})-(\ref{mKdV2}). Then $p(t,x):=\log
\psi(t,x)$ satisfies \eqref{EDE-asymptotics1}-\eqref{EDE-asymptotics2} with
$p_{0}(x)=\int_{0}^{x}r_{0}(s)\;ds$. As $p_{0}$ is in $\mathcal{O}^{*}%
_{\beta+1}(\mathbb{R})$ the solution $p(t,x)$ of
\eqref{EDE-asymptotics1}-\eqref{EDE-asymptotics2} belongs to $\mathcal{O}%
^{*}_{\beta+1}(\mathbb{R}\times\mathbb{R})$ (cf. Lemma
\ref{Lem:EDE-asymptoticsO} below). In particular $r(t,x)=p_{x}(t,x)\in
\mathcal{O}_{\beta}(\mathbb{R}\times\mathbb{R})$.

The uniqueness of the solution $r(t,x)$ constructed above follows from Theorem
\ref{Th:mKdV} $(b)$ and the uniqueness result for KdV in Theorem 1 in
\cite{BS2} (cf. Theorem \ref{Th:SB-uniqueness}, Appendix B). \hfill$\Box$
\vspace{3mm}

In the proof of Theorem \ref{Th:mKdV-O} we used the following analogue of
Lemma \ref{Lem:EDE-asymptotics}.

\begin{Lemma}
\hspace{-2mm}\textbf{.}\label{Lem:EDE-asymptoticsO} Let $\beta<1/2$ and
$\delta:=\max\{2\beta,\beta-1\}<1$. Consider the initial value problem
\eqref{EDE-asymptotics1}-\eqref{EDE-asymptotics2} where $q\in\mathcal{O}%
_{\delta}(\mathbb{R}\times\mathbb{R})$. Then, for any initial data $p_{0}%
\in\mathcal{O}^{*}_{\beta+1}(\mathbb{R})$, there exists a solution
$p\in\mathcal{O}^{*}_{\beta+1}(\mathbb{R}\times\mathbb{R})$ of
\eqref{EDE-asymptotics1}-\eqref{EDE-asymptotics2} which is unique in
$C^{\infty}(\mathbb{R}\times\mathbb{R})$.
\end{Lemma}

\noindent\emph{Proof of Lemma \ref{Lem:EDE-asymptoticsO}.} The lemma is proved
by the same arguments as the ones used in the proof of Lemma \ref{Lem:EDE-S}
(see also \cite{BS3}, Proposition 1). \hfill$\Box$ \vspace{3mm}

\noindent\emph{Proof of Theorem \ref{Th:mKdV-o}.} The proof is similar to the
proof of Theorem \ref{Th:mKdV-O} and is based on the existence and uniqueness
results for the initial value problem of KdV of \cite{BS2,BS3} (cf. Theorem
\ref{Th:SB-uniqueness}, \ref{Th:SB-o} in Appendix B) and on a variant of Lemma
\ref{Lem:EDE-asymptoticsO} where the spaces $\mathcal{O}_{\beta}$ and
$\mathcal{O}^{*}_{\beta+1}$ are replaced by the spaces $o_{\beta}$ and
$o^{*}_{\beta+1}$ respectively. \hfill$\Box$ \vspace{3mm}

We conclude this section by stating a more general uniqueness result for the
mKdV initial value problem (\ref{mKdV1})-(\ref{mKdV2}). Let $I=(a,b)\subseteq
\mathbb{R}$ with $-\infty\le a<b\le\infty$. Denote by $\mathcal{G}%
(\mathbb{R}\times\mathbb{R})$ the linear space of functions $r(t,x)$ in
$C^{\infty}(\mathbb{R}\times\mathbb{R})$ such that for any compact interval
$J\subseteq I$ one has for $|x|\ge1$ and any $k\ge1$
\[
r(t,x)=o(\sqrt{|x|})\;\;\;\mbox{and}\;\;\;\partial^{k}_{x}r(t,x)=O(1/\sqrt
{|x|}),
\]
uniformly in $t\in J$. The following theorem follows in a straightforward way
from Theorem \ref{Th:mKdV} $(b)$ and Theorem 1 in \cite{BS2} (cf. Theorem
\ref{Th:SB-uniqueness}, Appendix B).

\begin{Th}
\hspace{-2mm}\textbf{.}\label{Th:mKdV-uniqueness} There exists at most one
solution of the mKdV initial value problem (\ref{mKdV1})-(\ref{mKdV2}) in
$\mathcal{G}(\mathbb{R}\times\mathbb{R})$.
\end{Th}

\section{Spectral invariance}

\label{Sec:Spect_Invariance} In this section, we prove the spectral invariance
of the Schr{\"odinger} operator $L(t)=-\partial_{x}^{2}+q(t,x)$ under the KdV
flow (see Theorem \ref{Prop:spectral_invariance} stated in the introduction)
and the spectral invariance of the impedance operator $T(t)=-\frac{d^{2}%
}{dx^{2}}-2r(t,x)\,\frac{d}{dx}$ under the mKdV flow (Theorem
\ref{Prop:spectral_invariance*} below).

\vspace{0.3cm}

\noindent\emph{Schr{\"o}dinger operator:} We prove the two statements of
Theorem \ref{Prop:spectral_invariance} separately, using two different
methods. In the case $\beta<1$, statement $(i)$ follows from statement $(ii)$.
We point out that our proof of statement $(i)$ is elementary and
self-contained, whereas the proof of statement $(ii)$ relies on a result from
\cite{BS2}.

Let $q\in\mathcal{O}_{\beta}(I\times\mathbb{R})$ with $\beta\le1$ and
$I=(a,b)\subseteq\mathbb{R}$, $-\infty\le a< b\le+\infty$. Consider the
$t$-parameter family of self-adjoint operators in $L^{2}(\mathbb{R})$
\[
L(t):=-\partial_{x}^{2}+q(t,x)\,.
\]
As for any $t\in I$ there exists a constant $C=C(t)>0$ such that $|q(t,x)|\le
C\,|x|$ for $|x|\ge1$, by Sears' theorem (cf. \cite[Chapter II]{BeS}) the
symmetric operator $u\mapsto-u^{\prime\prime}+q(t) u$ in $L^{2}(\mathbb{R})$
with the domain $C^{\infty}_{0}(\mathbb{R})$ is essentially self-adjoint. We
define $L(t)$ to be the closure of this operator.

In order to prove Theorem \ref{Prop:spectral_invariance}, we first extend the
Q-L formalism to operators with spectral parameter. Let $q\in C^{\infty
}(I\times\mathbb{R})$. For a given $\lambda\in\mathbb{\mathbb{R}}$ consider
the operators
\[
Q_{\lambda}(t):=(4\lambda+2q(t,x))\partial_{x}-q_{x}(t,x)
\]
and
\[
L_{\lambda}(t):=L(t)-\lambda\,.
\]
The operators $Q_{\lambda}$ and $L_{\lambda}$ satisfy the commutator relation
\begin{equation}
\label{e:Q-A*}\dot{L}_{\lambda}=[Q_{\lambda},L_{\lambda}]+4q_{x} L_{\lambda
}+KdV(q)\,.
\end{equation}
In particular, $q(t,x)$ is a solution of the KdV equation, $KdV(q)=0$, if and
only if
\begin{equation}
\label{e:Q-A}\dot{L}_{\lambda}=[Q_{\lambda},L_{\lambda}]+4q_{x} L_{\lambda}%
\end{equation}
for some $\lambda\in\mathbb{R}$ (and hence, for all $\lambda\in\mathbb{R}$).
Now, assume that $q\in\mathcal{O}_{\beta}(I\times\mathbb{R})$ with $\beta\le1$
is a solution of the KdV equation and let $\psi\in C^{\infty}(I\times
\mathbb{R})$ be the solution of the first-order evolution equation
\begin{align}
\psi_{t}(t,x)  &  =(4\lambda+2q(t,x))\psi_{x}(t,x)-q_{x}(t,x)\psi
(t,x)\label{e:epde1*}\\
\psi(t^{\prime},x)  &  =\psi_{0}(x)\,. \label{e:epde2*}%
\end{align}
Let $\varphi(t,x):=L_{\lambda}(t)\psi(t,x)$. Using \eqref{e:Q-A} we get
\begin{align}
\varphi_{t}  &  ={\dot L}_{\lambda}\psi+L_{\lambda}\psi_{t}\nonumber\\
&  =([Q_{\lambda},L_{\lambda}]+4q_{x}L_{\lambda})\psi+ L_{\lambda}Q_{\lambda
}\psi\nonumber\\
&  =Q_{\lambda}(L_{\lambda}\psi)+4q_{x}(L_{\lambda}\psi)\nonumber\\
&  =(4\lambda+2q)\varphi_{x}+3q_{x}\varphi\;.\nonumber
\end{align}
Hence, $\varphi\equiv\varphi(t,x)$ is a solution of the initial value problem
(cf. Lemma \ref{Lem:EDE})
\begin{align}
\varphi_{t}(t,x)  &  =(4\lambda+2q(t,x))\varphi_{x}(t,x)+3q_{x}(t,x)\varphi
(t,x)\nonumber\\
\varphi|_{t=t^{\prime}}  &  =L_{\lambda}(t^{\prime})\psi_{0}\;.\nonumber
\end{align}

\begin{Lemma}
\hspace{-2mm}\textbf{.}\label{Lem:inequalities} Let $q\in\mathcal{O}_{\beta
}(I\times\mathbb{R})$ with $\beta\le1$ be a solution of the KdV equation and
let $J$ be a compact interval in $I$. Then there exists a constant $C(J)>0$
such that for any $t,t^{\prime}\in J$ and for any $\psi_{0}\in C^{\infty}%
_{0}(\mathbb{R})$ the solution $\psi$ of the evolution equation
\eqref{e:epde1*}-\eqref{e:epde2*} satisfies the inequalities
\begin{equation}
\label{e:inequality1}\|\psi(t)\|\le C(J)\,\|\psi_{0}\|
\end{equation}
and
\begin{equation}
\label{inequality2}\|L_{\lambda}(t)\psi(t)\|\le C(J)\,\|L_{\lambda}(t^{\prime
})\psi_{0}\|\,.
\end{equation}

\end{Lemma}

\noindent\emph{Proof of Lemma \ref{Lem:inequalities}.} Let $\xi(t;t^{\prime
},x_{0})$ be the solution of the ordinary differential equation
\begin{align}
\dot{\xi}=-(4\lambda+2q(t,\xi))\label{e:ode1}\\
\xi|_{t=t^{\prime}}=x \,. \label{e:ode2}%
\end{align}
It follows from the equation of variation of \eqref{e:ode1}-\eqref{e:ode2}
that
\begin{equation}
\label{e:xi_x}\xi_{x}(t;t^{\prime},x)=e^{-2\int_{t^{\prime}}^{t}q_{x}(\tau
,\xi(\tau;t^{\prime},x))\,d\tau}\,.
\end{equation}
By the method of characteristics we get that for any $t^{\prime\prime}\in I$,
\[
\phi(t^{\prime\prime},\xi(t^{\prime\prime};t^{\prime},x))=\phi_{0}%
(x)\,e^{3\int_{t^{\prime}}^{t^{\prime\prime}}q_{x}(\tau,\xi(\tau;t^{\prime
},x))\,d\tau}\,,
\]
or, equivalently,
\begin{equation}
\label{e:phi2}\phi(t^{\prime\prime},y)=\phi_{0}(\xi(t^{\prime};t^{\prime
\prime},y))\,e^{3\int_{t^{\prime}}^{t^{\prime\prime}}q_{x}(\tau,\xi
(\tau;t^{\prime\prime},y))\,d\tau}\,.
\end{equation}
Assuming that $t^{\prime},t^{\prime\prime}\in J$ for some compact interval $J$
in $I$, we obtain from \eqref{e:xi_x}, \eqref{e:phi2}, and $q\in
\mathcal{O}_{\beta}(I\times\mathbb{R})$ with $\beta\le1$ that there exist
constants $C(J)>0$ and $C_{1}(J)>0$ such that for any $\phi_{0}\in C^{\infty
}_{0}(\mathbb{R})$
\begin{align}
\|\phi(t^{\prime\prime})\|^{2}  &  =\int_{-\infty}^{\infty}\phi_{0}%
(\xi(t^{\prime};t^{\prime\prime},y))^{2}\, e^{6\int_{t^{\prime}}%
^{t^{\prime\prime}}q_{x}(\tau,\xi(\tau;t^{\prime\prime},y))\,d\tau
}\,dy\nonumber\\
&  \le C_{1}(J)\,\int_{-\infty}^{\infty}\phi_{0}(\xi(t^{\prime};t^{\prime
\prime},y))^{2}\,dy\nonumber\\
&  =C_{1}(J)\,\int_{-\infty}^{\infty}\phi_{0}(x)^{2}|\xi_{x}(t^{\prime\prime
};t^{\prime},x)|\,dx\nonumber\\
&  \le C(J)^{2}\,\|\phi_{0}\|^{2}\,.\nonumber
\end{align}
Arguing as above and choosing $C(J)>0$ larger if necessary, one proves that
for any $t^{\prime\prime}\in J$ and for any $\psi_{0}\in C^{\infty}%
_{0}(\mathbb{R})$ the solution $\psi(t,x)$ of
\eqref{e:epde1*}-\eqref{e:epde2*} satisfies
\[
\|\psi(t^{\prime\prime})\|\le C(J)\,\|\psi_{0}\|\,.
\]
This completes the proof of the lemma. \hfill$\Box$ \vspace{3mm}

\noindent\emph{Proof of statement $(i)$ of Theorem
\ref{Prop:spectral_invariance}.} A point $\lambda\in\mathbb{R}$ belongs to the
spectrum of the self-adjoint operator $L(t)$ iff there exist a sequence
$(\varepsilon_{k})_{k\ge1}$ of positive numbers $\varepsilon_{k}>0$ with
$\lim\limits_{k\to\infty}\varepsilon_{k}=0$ and a sequence of functions
$(\psi_{k})_{k\ge1}\subseteq C^{\infty}_{0}(\mathbb{R})$ , $\psi
_{k}\not \equiv 0$, such that for any $k\ge1$,
\begin{equation}
\label{e:lambda_in_spec}\|L_{\lambda}(t)\psi_{k}\|\le\varepsilon_{k}\|\psi
_{k}\|\,.
\end{equation}
Assume that $\lambda\in\mathop{\tt spec} L(t^{\prime})$. We will prove that
for any $t^{\prime\prime}\in I$, $\lambda\in\mathop{\tt spec} L(t^{\prime
\prime})$. Take $(\epsilon_{k})$ and $(\psi_{k})$ as above so that
\eqref{e:lambda_in_spec} is satisfied with $t=t^{\prime}$. Using Lemma
\ref{Lem:inequalities} we get the following estimates for the solution
$\psi_{k}(t,x)$ of equation \eqref{e:epde1*} with initial data $\psi
_{k}|_{t=t^{\prime}}=\psi_{k}$
\begin{align}
\|L(t^{\prime\prime})\psi_{k}(t^{\prime\prime})\|  &  \le C(J)\,\|L(t^{\prime
})\psi_{k}\|\nonumber\\
&  \le\epsilon_{k} C(J)\,\|\psi_{k}\|\le\epsilon_{k} C(J)^{2}\|\psi
_{k}(t^{\prime\prime})\|\,. \label{e:*}%
\end{align}
As $\epsilon_{k} C(J)^{2}\to0$ as $k\to\infty$ we get from \eqref{e:*} that
$\lambda\in\mathop{\tt spec} L(t^{\prime\prime})$. As the inclusion
$\mathop{\tt spec} L(t^{\prime})\subseteq\mathop{\tt spec} L(t^{\prime\prime
})$ was proved for any $t^{\prime},t^{\prime\prime}\in I$, $\mathop{\tt spec}
L(t^{\prime})=\mathop{\tt spec} L(t^{\prime\prime})$.

Any eigenfunction of $L(t^{\prime})$ with eigenvalue $\lambda\in\mathbb{R}$
coincides up to a set of measure zero with a smooth solution of the
differential equation $-\psi_{0}^{\prime\prime}(x)+q(t,x)\psi_{0}%
(x)=\lambda\psi_{0}(x)$ such that $\psi_{0}\in L^{2}(\mathbb{R})$. Let $\psi$
be the solution of the evolution equation \eqref{e:epde1*}-\eqref{e:epde2*}.
Arguing as in the proof of Lemma \ref{Lem:inequalities} one sees that the
inequalities \eqref{e:inequality1} and \eqref{inequality2} still hold. In
particular, we get that for any $t\in I$, $\psi(t)\in L^{2}(\mathbb{R})$ and
$L(t) \psi(t)=\lambda\psi(t)$. The coincidence of the multiplicities follows
from the uniqueness of solution for the initial value problem
\eqref{e:epde1*}-\eqref{e:epde2*}. \hfill$\Box$ \vspace{3mm}

\begin{Rem}
\hspace{-2mm}\textbf{.} In addition to the commutator relation \eqref{e:Q-A*},
the operators $L_{\lambda}$ and $Q_{\lambda}$ satisfy for arbitrary
$\lambda\in\mathbb{R}$ the identity
\[
L_{\lambda}(\partial_{t}-Q_{\lambda})+[L_{\lambda}(\partial_{t}-Q_{\lambda
})]^{*}=-KdV(q)
\]
where $P^{*}$ denotes the formal adjoint of a differential operator $P$ in
$L^{2}(I\times\mathbb{R})$.
\end{Rem}

\vspace{0.3cm}

\noindent\emph{Proof of statement $(ii)$ of Theorem
\ref{Prop:spectral_invariance}.} Assume either $q\in\mathcal{O}_{\beta
}(I\times\mathbb{R})$ with $\beta< 1$ or $q\in{o}_{\beta}(I\times\mathbb{R})$
with $\beta\le1$ where $I=(a,b)\subseteq\mathbb{R}$ and $-\infty\le a<
b\le+\infty$. In addition, assume for simplicity that $t^{\prime}=0$.
Following Lax \cite{Lax} we consider the one parameter family of third-order
linear differential operators
\begin{equation}
\label{e:A}A(t):=-4\partial_{x}^{3}+6q(t,x)\partial_{x}+3q_{x}(t,x)
\end{equation}
with $t$ as parameter. The operators $L(t)=-\partial_{x}^{2}+q(t,x)$ and
$A(t)$ satisfy the commutator relation
\[
{\dot L}=[A,L]+KdV(q)
\]
where $KdV(q)=q_{t}-6q q_{x}+q_{xxx}$. Assuming that $q(t,x)$ is a solution of
the KdV equation we obtain that $L(t)$ and $A(t)$ satisfy the classical Lax
pair relation
\begin{equation}
\label{e:Lax_pair}{\dot L}=[A,L]\,.
\end{equation}

Consider the linear evolution equation
\begin{align}
\psi_{t}(t)  &  =A(t)\psi(t)\label{e:Lax-evolution1}\\
\psi|_{t=0}  &  =\psi_{0} \label{e:Lax-evolution2}%
\end{align}
with initial data $\psi_{0}$ in the Schwartz space $\mathcal{S}(\mathbb{R})$.
The existence of a solution of
\eqref{e:Lax-evolution1}-\eqref{e:Lax-evolution2} evolving in $\mathcal{S}%
(\mathbb{R})$ on the whole time interval $(a,b)$ was solved positively in
\cite{BS2} by applying a difference scheme method as in \cite{B1,M}. According
to \cite[Theorem 2]{BS2}, \eqref{e:Lax-evolution1}-\eqref{e:Lax-evolution2}
has a solution $\psi\in C^{1}(I,\mathcal{S}(\mathbb{R}))$.\footnote{In fact,
this solution lies in $C^{\infty}(I,\mathcal{S}(\mathbb{R}))$.} By applying
Holmgren's principle one sees that the solution $\psi$ is indeed unique.
Denote by $\Psi(t)$ the operator
\begin{equation}
\label{e:Psi}\Psi(t) : \mathcal{S}(\mathbb{R})\to\mathcal{S}(\mathbb{R}%
),\;\;\;\psi_{0}\mapsto\psi(t).
\end{equation}
Using that $A(t)$ is skew-symmetric for any $t$ one gets by integration by
parts, that for any $\psi_{0}\in\mathcal{S}(\mathbb{R})$
\[
\frac{d}{dt}\Big(\Psi(t)\psi_{0},\Psi(t)\psi_{0}\Big)=\Big(A(t)\Psi(t)\psi
_{0},\Psi(t)\psi_{0}\Big)+\Big(\Psi(t)\psi_{0},A(t)\Psi(t)\psi_{0}\Big)=0
\]
where $(\cdot,\cdot)$ denotes the $L^{2}$-scalar product. Hence, the operator
\eqref{e:Psi} preserves the $L^{2}$-norm. As $\mathcal{S}(\mathbb{R})$ is
dense in $L^{2}(\mathbb{R})$, $\Psi(t)$ then extends to a unitary operator on
$L^{2}(\mathbb{R})$.

Let $\psi(t)$ be the solution of
\eqref{e:Lax-evolution1}-\eqref{e:Lax-evolution2}. By \eqref{e:Lax_pair} and
the Leibniz rule we get
\begin{align}
(L\psi)^{\cdot}={\dot L}\psi+L{\dot\psi}=[A,L]\psi+LA\psi=A(L\psi)\,.
\end{align}
Hence, $L(t)\psi(t)$ is a solution of \eqref{e:Lax-evolution1} with initial
data $L(0)\psi_{0}$. The latter result together with the uniqueness of the
solution of \eqref{e:Lax-evolution1} with the initial data $L(0)\psi_{0}$
imply that
\begin{equation}
\label{e:conjugation}\Psi(t)L(0)\psi_{0}=L(t)\Psi(t)\psi_{0}\;\;\;\forall
\;\psi_{0}\in\mathcal{S}(\mathbb{R})\,.
\end{equation}
Since, by Sears' theorem, $\mathcal{S}(\mathbb{R})$ is dense in the domains of
$L(0)$ and $L(t)$ with respect to the graph norms, it follows from
\eqref{e:Psi} that the identity \eqref{e:conjugation} holds for any $\psi_{0}$
in the domain of $L(0)$. Together with the property that $\Psi(t)$ is unitary,
this establishes the claimed unitary equivalence of the operators $L(0)$ and
$L(t)$. \hfill$\Box$ \vspace{3mm}

\vspace{0.3cm}

\noindent\emph{Impedance operator:} Here we prove that the spectrum of the
impedance operator
\[
T(t):=-\frac{d^{2}}{dx^{2}}-2r(t,x)\,\frac{d}{dx}
\]
where $r(t,x)$ is a solution of the mKdV equation in $\mathcal{O}_{\beta
}(I\times\mathbb{R})$ with $\beta\le1/2$ and $I=(a,b)$, $-\infty\le
a<b\le+\infty$, is invariant. For any given $t\in(a,b)$ the operator
$-\frac{d^{2}}{dx^{2}}-2r(t,x)\,\frac{d}{dx}$ with domain $C^{\infty}%
_{0}(\mathbb{R})$ is essentially self-adjoint in the weighted $L^{2}$-space
$L^{2}(\mathbb{R},\rho^{2}dx)$ with the density function $\rho(t,x)^{2}%
:=e^{2\int_{0}^{x}r(t,s)\,ds}$ (see below). We denote by $T(t)$ the closure of
this operator in $L^{2}(\mathbb{R},\rho^{2}dx)$.

\begin{Th}
\hspace{-2mm}\textbf{.}\label{Prop:spectral_invariance*} Let $r\in
\mathcal{O}_{\beta}(I\times\mathbb{R})$ with $\beta\le1/2$ be a solution of
the mKdV equation. Then for any $t,t^{\prime}\in I$
\[
\mathop{\tt spec} T(t)= \mathop{\tt spec} T(t^{\prime})\,.
\]
Moreover, if $\beta<1/2$ or $r\in o_{\beta}(I\times\mathbb{R})$ with $\beta
\le1/2$, then the operators $T(t)$ and $T(t^{\prime})$ are unitarily equivalent.
\end{Th}

To prove Theorem \ref{Prop:spectral_invariance*} we first need to establish
the following auxiliary result.

\begin{Lemma}
\hspace{-2mm}\textbf{.}\label{Lem:diagram} Let $r\in C^{\infty}(\mathbb{R})$,
$q(x)=r^{\prime}+r^{2}$, and $\rho(x)=e^{\int_{0}^{x}r(s)\,ds}$.

\begin{itemize}
\item[(a)] The map
\[
\Phi_{\rho}: C^{\infty}_{0}(\mathbb{R})\to C^{\infty}_{0}(\mathbb{R}%
),\;\;\;u(x)\mapsto\rho(x)\,u(x)
\]
extends to an isometry $L^{2}(\mathbb{R},\rho^{2}dx)\to L^{2}(\mathbb{R})$.

\item[(b)] The diagram
\[%
\begin{array}
[c]{rcccccl}%
L^{2}(\mathbb{R},\rho^{2}dx) & \supseteq & C^{\infty}_{0}(\mathbb{R}) &
\overset{T_{r}}{\longrightarrow} & C^{\infty}_{0}(\mathbb{R}) & \subseteq &
L^{2}(\mathbb{R},\rho^{2}dx)\\
&  & \downarrow\lefteqn{\Phi_{\rho}} &  & \downarrow\lefteqn{\Phi_{\rho}} &  &
\\
L^{2}(\mathbb{R}) & \supseteq & C^{\infty}_{0}(\mathbb{R}) & \overset{L_{q}%
}{\longrightarrow} & C^{\infty}_{0}(\mathbb{R}) & \subseteq & L^{2}%
(\mathbb{R})\\
&  &  &  &  &  &
\end{array}
\]
where $T_{r}=-\frac{d^{2}}{dx^{2}}-2r\,\frac{d}{dx}$ and $L_{q}=-\frac{d^{2}%
}{dx^{2}}+q$, is commutative.
\end{itemize}
\end{Lemma}

\noindent\emph{Proof of Lemma \ref{Lem:diagram}.} $(a)$ is obvious. To prove
$(b)$ use the relation $\rho^{\prime}/\rho=r$ to get for any $u\in C^{\infty
}_{0}(\mathbb{R})$
\[
T_{r}(u)=-(\rho^{2}\,u^{\prime})^{\prime}/\rho^{2}\,.
\]
Hence, for any $w\in C^{\infty}_{0}(\mathbb{R})$,
\begin{align}
T_{r}\circ\Phi_{\rho}^{-1}(w)  &  =T_{r}(w/\rho)\nonumber\\
&  =-(w^{\prime\prime}\rho-w\rho^{\prime\prime})/\rho^{2}\nonumber\\
&  =(-w^{\prime\prime}+(\rho^{\prime\prime}/\rho)\,w)/\rho\nonumber\\
&  =\Phi_{\rho}^{-1}\circ L_{q}(w)\nonumber
\end{align}
where we have used that $\frac{\rho^{\prime\prime}}{\rho}=\Big(\frac
{\rho^{\prime}}{\rho}\Big)^{\prime}+\Big(\frac{\rho^{\prime}}{\rho}\Big)^{2}$.
\hfill$\Box$ \vspace{3mm}

\noindent Assume that $r\in\mathcal{O}_{\beta}(\mathbb{R})$ with $\beta\le
1/2$. Then $q=r^{\prime}+r^{2}\in\mathcal{O}_{\delta}(\mathbb{R})$ with
$\delta\le1$. As, by Sears' theorem, the Schr{\"o}dinger operator $L_{q}$ is
essentially self-adjoint, the impedance operator $T_{r}$ is essentially
self-adjoint by Lemma \ref{Lem:diagram}. Moreover, it follows from Lemma
\ref{Lem:diagram} that the closures of both operators are unitarily
equivalent. In particular,
\begin{equation}
\label{e:isospectral}\mathop{\tt spec} T_{r}=\mathop{\tt spec} L_{q}\,.
\end{equation}
\noindent\emph{Proof of Theorem \ref{Prop:spectral_invariance*}.} The
statement of the theorem follows from the unitary equivalence of $L_{q}$ and
$T_{r}$, Theorem \ref{Prop:spectral_invariance}, and the fact that the Miura
map $r\mapsto r_{x}+r^{2}$ maps smooth solutions of mKdV to smooth solutions
of KdV. \hfill$\Box$ \vspace{3mm}

\section{Evolution of generalized eigenfunctions}

In this section, we consider the family of first-order differential operators

\noindent$Q_{\lambda}(t)=(4\lambda+2q(t,x))\partial_{x}-q_{x}(t,x)$, where $t$
is a real parameter in the interval $I=(a,b)$, $\lambda$ is the spectral
parameter and $q(t,x)$ is a solution of the KdV equation in $\mathcal{O}%
_{\beta}(I\times\mathbb{R})$ with $\beta<1$ or in $o_{\beta}(I\times
\mathbb{R})$ with $\beta\le1$. We then prove that for any $t^{\prime}$,
$t^{\prime\prime}$ in $I$ the solution operator corresponding to the family of
evolution equations induced by $Q_{\lambda}(t)$, transforms any complete
orthonormal system of generalized eigenfunctions of the operator $L(t^{\prime
})$ to such a system for the operator $L(t^{\prime\prime})$ (Theorem
\ref{Th:generalized_eigenfunctions}). Moreover, in Proposition
\ref{Prop:evolution} below, we prove that the solution of the evolution
equation involving the third-order differential operator $A(t)$ appearing in
the classical Lax pair for KdV (see \eqref{e:A}), can be obtained in terms of
the solutions of the one-parameter family of first-order evolution equations
$\psi_{t}=Q_{\lambda}(t)\psi$ where $\lambda\in\mathbb{R}$ is a spectral parameter.

Assume that $q$ is a solution of KdV satisfying either $q\in\mathcal{O}%
_{\beta}(I\times\mathbb{R})$ with $\beta< 1$ or $q\in{o}_{\beta}%
(I\times\mathbb{R})$ with $\beta\le1$ where $I=(a,b)$, $-\infty\le
a<b\le+\infty$ and for convenience, $0\in(a,b)$. Let us consider the family of
the first-order evolution equations
\begin{align}
\psi_{t}(t)  &  =Q_{\lambda}(t)\psi(t)\label{LEE1}\\
\psi|_{t=0}  &  =\psi_{0} \label{LEE2}%
\end{align}
where $\lambda\in\mathbb{R}$ is a parameter. According to Lemma \ref{Lem:EDE},
for any $\psi_{0}\in C^{\infty}(\mathbb{R})$, the initial value problem
\eqref{LEE1}-\eqref{LEE2} has a unique solution in $C^{\infty}(I\times
\mathbb{R})$. Moreover, if the initial data $\psi_{0}$ lies in the Schwartz
space $\mathcal{S}(\mathbb{R})$ then $\psi(t)\in\mathcal{S}(\mathbb{R})$ for
any $t\in I$, and $\psi\in C^{1}(I,\mathcal{S}(\mathbb{R}))$ (see Lemma
\ref{Lem:EDE-S}).

Consider a Hilbert-Schmidt rigging
\[
{\mathcal{H}}_{+}\subseteq{\mathcal{H}}=L^{2}(\mathbb{R})\subseteq
{\mathcal{H}}_{-},
\]
associated to the Hilbert-Schmidt operator
\[
K: L^{2}(\mathbb{R})\to L^{2}(\mathbb{R}),\;\;\;K:=(-\partial_{x}^{2}%
+x^{2})^{-s}\,,
\]
with some $s>1/2$. Recall that ${\mathcal{H}}_{+}:=K({\mathcal{H}})$ and the
norm $\|\cdot\|_{+}$ in ${\mathcal{H}}_{+}$ is defined by
\[
\|Kh\|_{+}=\|h\|\,\;\;\forall h\in{\mathcal{H}},
\]
whereas ${\mathcal{H}}_{-}$ is the dual to ${\mathcal{H}}_{+}$ (see
\cite[Supplement 1, \S \,3, 4]{BeS}). One has the following chain of
continuous embeddings
\begin{equation}
\label{e:rigging}\mathcal{S}(\mathbb{R})\subseteq{\mathcal{H}}_{+}\subseteq
L^{2}(\mathbb{R})\subseteq{\mathcal{H}}_{-}\subseteq\mathcal{S}^{\prime
}(\mathbb{R})
\end{equation}
where $\mathcal{S}^{\prime}(\mathbb{R})$ denotes the space of tempered
distributions. It follows from \cite[Supplement 1, \S \,7]{BeS} that there
exists a complete orthonormal system of generalized eigenfunctions
\begin{equation}
\label{e:generalized_eigenfunctions_t=0}\{\psi(x,m)\,|\,m\in M\}\subseteq
C^{\infty}(\mathbb{R})\cap{\mathcal{H}}_{-}\subseteq\mathcal{S}^{\prime
}(\mathbb{R})\,,
\end{equation}
of $L(0)$ with generalized eigenvalues $\lambda(m)$, $m\in M$, on a measure
space $M$ with the measure $\mu$. By definition, the functions
\eqref{e:generalized_eigenfunctions_t=0} constitute a complete orthonormal
system of generalized eigenfunctions if the following two properties hold (cf.
\cite[Supplement 1, Definition 2.4]{BeS}):

\begin{itemize}
\item[(i)] for any $h_{+}\in{\mathcal{H}}_{+}$, the function $m\mapsto
(h_{+},\psi(m))$ on $M$ belongs to $L^{2}(M,\mu)$;

\item[(ii)] the map $h_{+}\mapsto F(h_{+})$, $F(h_{+})(m):=(h_{+},\psi(m))$,
extends to a unitary operator ${\mathcal{H}}\to L^{2}(M,\mu)$.
\end{itemize}

Here and below we use the notation $(\cdot,\cdot)$ for miscellaneous
sesquilinear dualities which extend the usual $L^{2}$-inner product by
continuity. So these dualities are linear with respect to the first argument
and antilinear with respect to the second one. The unitary transform $F :
{\mathcal{H}}\to L^{2}(M,\mu)$ is called \emph{generalized Fourier transform}
corresponding to the system \eqref{e:generalized_eigenfunctions_t=0}. Note
that the functions $\psi(x,m)$ are smooth in $x$ and satisfy the relation
\begin{equation}
\label{e:eigenfunction}L(0)\psi(m)=\lambda(m)\psi(m)
\end{equation}
where $L(0)$ is the differential operator $-\partial_{x}^{2}+q(0,x)$.

\begin{Rem}
\hspace{-2mm}\textbf{.} It follows from \eqref{e:eigenfunction} and the
construction of generalized eigenfunctions in \cite[Supplement 1,\S \,5]{BeS}
that the multiplicity of the spectrum of $L(0)$ is at most two. Hence we can
choose $M$ to be the disjoint union of two copies of $\mathbb{R}$,
$M=\mathbb{R}\sqcup\mathbb{R}$, with positive finite Lebesgue-Stiltjes measure
on each of the two copies.
\end{Rem}

\begin{Th}
\hspace{-2mm}\textbf{.}\label{Th:generalized_eigenfunctions} Denote by
$\psi(t,x,m)$ the solution of the initial value problem
\eqref{LEE1}-\eqref{LEE2} with the initial data $\psi_{0}(x)=\psi(x,m)$ and
$\lambda=\lambda(m)$. Then for any $t\in I$, $\{\psi(t,x,m)\,|\,m\in M\}$ is a
complete orthonormal system of generalized eigenfunctions of the operator
$L(t)$ on the measure space $M$ with the same measure $\mu$ and the same
generalized eigenvalues $\lambda(m)$.
\end{Th}

To prove Theorem \ref{Th:generalized_eigenfunctions} we first need to consider
the linear evolution equation in $\mathcal{S}^{\prime}(\mathbb{R})$
\begin{align}
\psi_{t}(t)  &  ={\hat A}(t)\psi(t)\label{e:Lax-evolution1'}\\
\psi|_{t=0}  &  =\psi_{0}\in\mathcal{S}^{\prime}(\mathbb{R})
\label{e:Lax-evolution2'}%
\end{align}
where ${\hat A}(t) : \mathcal{S}^{\prime}(\mathbb{R})\to\mathcal{S}^{\prime
}(\mathbb{R})$ denotes the extension by continuity of the operator \eqref{e:A}
from $\mathcal{S}(\mathbb{R})$ to $\mathcal{S}^{\prime}(\mathbb{R})$ i.e.,
$\forall\;\psi\in\mathcal{S}^{\prime}(\mathbb{R})$ and $\forall\varphi
\in\mathcal{S}(\mathbb{R})$, $({\hat A}(t)\psi,\varphi):=(\psi,A^{*}%
(t)\varphi)$ where $A^{*}(t)$ is the formal adjoint to $A(t)$, which coincides
with $-A(t)$, and the derivative $\psi_{t}=d\psi/dt$ is understood in the weak
topology of $\mathcal{S}^{\prime}(\mathbb{R})$.

\begin{Lemma}
\hspace{-2mm}\textbf{.}\label{Lem:weak_evolution} The initial value problem
\eqref{e:Lax-evolution1'}-\eqref{e:Lax-evolution2'} has unique solution in
$\mathcal{S}^{\prime}(\mathbb{R})$.
\end{Lemma}

\noindent\emph{Proof of Lemma \ref{Lem:weak_evolution}.} Let us first prove
the existence. Consider again the operator \eqref{e:Psi} extended by
continuity to an isometry in $L^{2}(\mathbb{R})$. Then define the curve
\begin{equation}
\label{e:solution}\psi: I\to\mathcal{S}^{\prime}(\mathbb{R}),\;t\mapsto
{\hat\Psi}(t)\psi_{0}\in\mathcal{S}^{\prime}(\mathbb{R})
\end{equation}
where $({\hat\Psi}(t)\chi,\varphi):=(\chi,\Psi(t)^{*}\varphi)$, $\forall
\,\chi\in\mathcal{S}^{\prime}(\mathbb{R})$, $\forall\,\varphi\in
\mathcal{S}(\mathbb{R})$. The operator $\Psi^{*}(t) : L^{2}(\mathbb{R})\to
L^{2}(\mathbb{R})$ denotes the adjoint operator of $\Psi(t) : L^{2}%
(\mathbb{R})\to L^{2}(\mathbb{R})$ with respect to the $L^{2}$-scalar product
and, as $\Psi(t)$ is unitary, coincides with $\Psi(t)^{-1}$. For any test
function $\varphi\in\mathcal{S}(\mathbb{R})$ one has
\[
\frac{d}{dt}(\psi(t),\varphi)= \frac{d}{dt}(\psi_{0},\Psi^{*}(t)\varphi)=
(\psi_{0},\Psi^{*}(t)A^{*}(t)\varphi)= ({\hat A}(t)\psi(t),\varphi)\,.
\]
Hence, $\psi(t)$ is a solution of \eqref{e:Lax-evolution1'}-\eqref{e:Lax-evolution2'}.

The uniqueness follows from Holmgren's principle and the existence of a
solution $\psi\in C^{1}(I,\mathcal{S}(\mathbb{R}))$ of the initial value
problem
\begin{align*}
\psi_{t}(t)  &  =-A(t)^{*}\psi(t)=A(t)\psi(t)\\
\psi|_{t=t_{0}}  &  =\psi_{0}%
\end{align*}
for any $t_{0}\in I$ (\cite[Theorem 2]{BS2}). \hfill$\Box$ \vspace{3mm}

In addition we need to consider the initial value problem in $\mathcal{S}%
^{\prime}(\mathbb{R})$
\begin{align}
\psi_{t}(t)  &  ={\hat Q}_{\lambda}(t)\psi(t)\label{e:Q-evolution1}\\
\psi|_{t=0}  &  =\psi_{0}\in\mathcal{S}^{\prime}(\mathbb{R}) \,,
\label{e:Q-evolution2}%
\end{align}
where ${\hat Q}_{\lambda}(t)$ is the extension by continuity of the
differential operator $Q_{\lambda}(t) : \mathcal{S}(\mathbb{R})\to
\mathcal{S}(\mathbb{R})$ to an operator on $\mathcal{S}^{\prime}(\mathbb{R})$.
Arguing as in the proof of Lemma \ref{Lem:weak_evolution}, one shows the
following lemma.

\begin{Lemma}
\hspace{-2mm}\textbf{.}\label{Lem:weak_evolution*} The initial value problem
\eqref{e:Q-evolution1}-\eqref{e:Q-evolution2} has unique solution in
$\mathcal{S}^{\prime}(\mathbb{R})$.
\end{Lemma}

Finally, for the proof of Theorem \ref{Th:generalized_eigenfunctions} we will
need the following result on unitary equivalence. Let $L_{1}$ and $L_{2}$ be
self-adjoint operators on a Hilbert space ${\mathcal{H}}$. Let ${\mathcal{H}%
}_{+}\subseteq{\mathcal{H}}\subseteq{\mathcal{H}}_{-}$ be a rigging associated
to a Hilbert-Schmidt operator $K : {\mathcal{H}}\to{\mathcal{H}}$ and let
\[
\{\psi_{1}(m)\,|\,m\in M\}\subseteq{\mathcal{H}}_{-}
\]
be a complete orthonormal system of generalized eigenfunctions of the operator
$L_{1}$ with generalized eigenvalues $\lambda(m)$, $m\in M$, on a measure
space $(M,\mu)$.

\begin{Lemma}
\hspace{-2mm}\textbf{.}\label{Lem:shift} Assume that the operators $L_{1}$ and
$L_{2}$ are unitarily equivalent via the isometry $\Psi: {\mathcal{H}}%
\to{\mathcal{H}}$, and define the system $\{\psi_{2}(m):={\hat\Psi}(\psi
_{1}(m))\,|\,m\in M\}$, where $({\hat\Psi}\chi,\phi):=(\chi,\Psi^{-1}\phi)$,
$\forall\chi\in{\mathcal{H}}_{-}$, $\forall\phi\in\Psi({\mathcal{H}}_{+})$.
Then $\{\psi_{2}(m)\,|\,m\in M \}$ is a complete orthonormal system of
generalized eigenfunctions of the operator $L_{2}$ on the same measure space
$(M,\mu)$ and with the same generalized eigenvalues $\lambda(m)$, $m\in M$
(with respect to the rigging $\Psi({\mathcal{H}}_{+})\subseteq{\mathcal{H}%
}\subseteq{\hat\Psi}({\mathcal{H}}_{-})$). Moreover,
\begin{equation}
\label{e:shift_relation}F_{2}(\phi)=F_{1}\circ\Psi^{-1}(\phi)
\end{equation}
where $F_{1}$ and $F_{2}$ denote the generalized Fourier transform
corresponding to the system $\{\psi_{1}(m)\,|\,m\in M\}$ and $\{\psi
_{2}(m)\,|\,m\in M\}$ respectively.
\end{Lemma}

\noindent\emph{Proof of Lemma \ref{Lem:shift}.} The proof of the lemma is
straightforward. Indeed, for any $\varphi\in{\mathcal{H}}_{+}$ one has a.e. on
$M$,
\[
F_{2}(\Psi(\varphi))(m)=(\Psi(\varphi),\psi_{2}(m))=(\Psi(\varphi),{\hat\Psi
}(\psi_{1}(m)))= (\varphi,\psi_{1}(m))=F_{1}(\varphi)(m)\,.
\]
In particular, one gets that for any $\phi\in\Psi({\mathcal{H}}_{+})$,
\eqref{e:shift_relation} holds. The latter relation proves that the
generalized Fourier transform $F_{2} : \Psi({\mathcal{H}}_{+})\to L^{2}%
(M,\mu)$ can be extended by continuity from $\Psi({\mathcal{H}}_{+})$ to an
isometry $F_{2} : {\mathcal{H}}\to L^{2}(M,\mu)$ that satisfies
\eqref{e:shift_relation} for any $\phi\in{\mathcal{H}}$. In particular, we get
that
\[
F_{2}^{-1}\circ{\hat\lambda}\circ F_{2}=\Psi\circ(F_{1}^{-1}\circ{\hat\lambda
}\circ F_{1})\circ\Psi^{-1}=\Psi\circ L_{1}\circ\Psi^{-1}=L_{1}
\]
where ${\hat\lambda}$ is the multiplication operator by $\lambda(m)$ in
$L^{2}(M,\mu)$. This completes the proof of the lemma. \hfill$\Box$
\vspace{3mm}

\noindent Now we are ready to prove Theorem
\ref{Th:generalized_eigenfunctions}.

\noindent\emph{Proof of Theorem \ref{Th:generalized_eigenfunctions}.} Assume
that $\psi(t,x,m)$ is the solution of the initial value problem
\eqref{LEE1}-\eqref{LEE2} with the initial data $\psi(x,m)$ and $\lambda
=\lambda(m)$ where $\psi(x,m)\in\mathcal{S}^{\prime}(\mathbb{R})\cap
C^{\infty}(\mathbb{R})$ is a generalized eigenfunction of $L(0)$ with the
generalized eigenvalue $\lambda(m)$. Then $I\to\mathcal{S}^{\prime}%
(\mathbb{R})$, $t\mapsto\psi(t,m)$, is the solution of \eqref{e:Q-evolution1}
with initial data $\psi(m)$. Using the commutator relation \eqref{e:Q-A*} with
$\lambda=\lambda(m)$ and arguing as in the proof of Proposition
\ref{Prop:eigenfunction-transport} one shows that $L_{\lambda(m)}%
(t)\psi(t,m)=0$. The latter together with the relation $A(t)=Q_{\lambda
}(t)+4\partial_{x}\circ L_{\lambda}(t)$ applied for $\lambda=\lambda(m)$,
implies
\begin{align*}
\psi(t,m)_{t}  &  ={\hat Q}_{\lambda(m)}(t)\psi(t,m)=({\hat A}(t)+4{\hat
\partial_{x}}\circ{\hat L}_{\lambda(m)}(t))\psi(t,m)\\
&  ={\hat A}(t)\psi(t,m)
\end{align*}
where ${\hat\partial_{x}}$ and ${\hat L}_{\lambda(m)}(t)$ are the extensions
by continuity of the differential operators $\partial_{x}$ and $L_{\lambda
(m)}(t)$ respectively. Hence, $t\mapsto\psi(t,m)$ solves
\eqref{e:Lax-evolution1'}-\eqref{e:Lax-evolution2'} with initial data
$\psi_{0}=\psi(m)$. By Lemma \ref{Lem:weak_evolution} and \eqref{e:solution},
we get that
\begin{equation}
\label{e:shift}\psi(t,m)={\hat\Psi}(t)(\psi(m))\,.
\end{equation}
By the proof of statement $(ii)$ of Theorem \ref{Prop:spectral_invariance}
(see Section \ref{Sec:Spect_Invariance}), $\Psi(t) : L^{2}(\mathbb{R})\to
L^{2}(\mathbb{R})$ is a unitary equivalence of the operators $L(0)$ and
$L(t)$. One then concludes from \eqref{e:shift} and Lemma \ref{Lem:shift} that
$\{\psi(t,m)\,|\,m\in M\}$ is a complete orthonormal system of generalized
eigenfunctions of $L(t)$ corresponding to the rigging of $L^{2}(\mathbb{R})$
obtained by shifting the rigging \eqref{e:rigging} via the isometry $\Psi(t) :
L^{2}(\mathbb{R})\to L^{2}(\mathbb{R})$. \hfill$\Box$ \vspace{3mm}

Denote by $F(0)$ the generalized Fourier transform corresponding to the
complete orthonormal system \eqref{e:generalized_eigenfunctions_t=0}, so that
for any initial data $\psi_{0} \in\mathcal{S}(\mathbb{R})$,
\begin{equation}
\label{E:Fourier}\tilde{\psi}_{0}(m) := (F(0) \psi_{0})(m) = \int_{\mathbb{R}}
\overline{\psi(x,m)} \psi_{0}(x) \, dx.
\end{equation}
Let $F(t)$ be the generalized Fourier transform corresponding to the system
$\{\psi(t,m)\,|\,m\in M\}$.

With the notation of Theorem \ref{Th:generalized_eigenfunctions}, we have the
following Corollary.

\begin{Coro}
\hspace{-2mm}\textbf{.}\label{Prop:evolution} For any $\psi_{0} \in
\mathcal{S}(\mathbb{R})$, the function
\[
\psi(t,x;\psi_{0}):=(F(t)^{-1}{\tilde\psi}_{0})(x)
\]
is a solution of the initial value problem \eqref{e:Lax-evolution1}-\eqref{e:Lax-evolution2}.
\end{Coro}

\begin{Rem}
\hspace{-2mm}\textbf{.} The solution obtained in Corollary
\ref{Prop:evolution} is unique in view of the discussion in the paragraph
following \eqref{e:Lax-evolution1}-\eqref{e:Lax-evolution2}.
\end{Rem}

\noindent\emph{Proof of Corollary \ref{Prop:evolution}.} As $\Psi(t) :
L^{2}(\mathbb{R})\to L^{2}(\mathbb{R})$ is a unitary equivalence of the
operators $L(0)$ and $L(t)$, it follows from \eqref{e:shift_relation} that
$F(0)=F(t)\circ\Psi(t)$. The latter relation leads to
\[
\psi(t;\psi_{0})=F(t)^{-1}\circ F(0)\psi_{0}=F(t)^{-1}\circ(F(t)\circ
\Psi(t))\psi_{0}=\Psi(t)\psi_{0}\,.
\]
As $\Psi(t)$ is the evolution operator of \eqref{e:Lax-evolution1},
$\psi(t;\psi_{0})$ solves \eqref{e:Lax-evolution1}-\eqref{e:Lax-evolution2}.
\hfill$\Box$ \vspace{3mm}

\begin{Rem}
\hspace{-2mm}\textbf{.}\label{Rem:FI} It follows from the construction of
generalized eigenfunctions of a self-adjoint operator in \cite[Supplement 1,
Theorem 2.1]{BeS} that there exists an isomorphism $U(t) : L^{2}%
(\mathbb{R})\to L^{2}(M,\mu)$ that transforms the system of generalized
eigenfunctions $\{\psi(t,x,m)\,|\,m\in M\}$ into a system of delta functions
$\{\chi(y,m):=\delta_{m}(y)\,|\,m\in M\}$ on $M$. (The delta function
$\delta_{m}$ is defined for almost every $m\in M$, and $\delta_{m}%
\in(F(0)({\mathcal{H}}_{+}))^{\prime}$ (see below).) The generalized Fourier
transform ${\tilde\psi}_{0}$ of the initial data $\psi_{0}$ belongs to the
space $F(0)({\mathcal{H}}_{+})$ that corresponds to the Hilbert-Schmidt
rigging $F(0)({\mathcal{H}}_{+})\subseteq L^{2}(M,\mu)\subseteq
(F(0)({\mathcal{H}}_{+}))^{\prime}$. Arguing as in the proof of
\cite[Supplement 1, Theorem 2.1]{BeS}, one sees that the system $\{\chi
(y,m)\,|\,m\in M\}$ is contained in the dual space to $F(0)({\mathcal{H}}%
_{+})$. Moreover, a trivial computation in the `model space' $L^{2}(M,\mu)$
leads to the relations \eqref{E:Fourier} and
\[
(F(t)^{-1}{\tilde\psi}_{0})(x)=\int_{M}\psi(t,x,m)\,{\tilde\psi}_{0}%
(m)\,d\mu(m),
\]
where the integral makes sense as a continuous extension of the $L^{2}$ scalar
product on $(M,\mu)$.
\end{Rem}

\begin{Rem}
\hspace{-2mm}\textbf{.} Corollary \ref{Prop:evolution} shows that
the solution of the evolution equation \eqref{e:Lax-evolution1} involving the
third-order differential operator $A(t)$ (see \eqref{e:A}), can be obtained in
terms of the solution of the evolution equation \eqref{LEE1}-\eqref{LEE2}
involving the first-order differential operator $Q_{\lambda}(t)=(4\lambda
+2q(t,x))\partial_{x}-q_{x}(t,x)$ with parameter $\lambda\in\mathbb{R}$,
provided we know the spectral decomposition of the operator $L(0)$.
\end{Rem}

\begin{Rem}
\hspace{-2mm}\textbf{.} There is another way to understand what
$\psi(x,m)$ is. We will assume that $M=\mathbb{R}\sqcup\mathbb{R}$, a disjoint
union of two lines $\mathbb{R}$, is a measure space with positive finite
Stiltjes measures on the lines (which is enough for us). Then the Schwartz
space $\mathcal{S}(\mathbb{R}\times M)$ is well defined, and so is
$\mathcal{S}^{\prime}$. Also the Schwartz kernel theorem works as well, and it
implies that every continuous linear operator from $\mathcal{S}(\mathbb{R})$
to $\mathcal{S}^{\prime}(M)$ or from $\mathcal{S}(M)$ to $\mathcal{S}^{\prime
}(\mathbb{R})$ (where $\mathcal{S}^{\prime}$ in both cases should be
considered with the weak topology), can be uniquely presented by a Schwartz
kernel from $\mathcal{S}^{\prime}(M\times\mathbb{R})$ or $\mathcal{S}^{\prime
}(\mathbb{R}\times M)$ respectievely. In particular this is true for all
bounded linear operators from $L^{2}(\mathbb{R})$ to $L^{2}(M)$ or back.

Vice versa, any distribution from $\mathcal{S}^{\prime}(M\times\mathbb{R})$
(or $\mathcal{S}^{\prime}(\mathbb{R}\times M)$) naturally defines a continuous
linear operator from $\mathcal{S}(\mathbb{R})$ to $\mathcal{S}^{\prime}(M)$
(or $\mathcal{S}(M)$ to $\mathcal{S}^{\prime}(\mathbb{R})$ respectively).

Now it is easy to see that $\psi=\psi(x,m)$ is a tempered distribution which
is the Schwartz kernel for the unitary operator $F$ (the generalized Fourier
transform) from $L^{2}(\mathbb{R})$ to $L^{2}(M)$, and $\overline{\psi(x,m)}$
is the Schwartz kernel for $F^{-1}$.
\end{Rem}

\appendix

\section{Appendix: Global solutions of a linear first-order PDE}

In this appendix we state and prove, for the convenience of the reader, a
result on the first-order linear PDE, used in the main body of the paper,
\begin{align}
u_{t}(t,x)  &  =a(t,x)u_{x}(t,x)+b(t,x)u(t,x)\label{EDE}\\
u|_{t=0}  &  =\psi(x) \label{EDE0}%
\end{align}
where $\psi\in C^{\infty}(\mathbb{R}) $, $a,b\in C^{\infty}(\mathbb{R}%
\times\mathbb{R})$, and $a$ grows for $x \to\pm\infty$ at most linearly. In
addition, we prove three technical lemmas used in the proof of Theorem
\ref{Th:mKdV-S}.

\begin{Lemma}
\hspace{-2mm}\textbf{.}\label{Lem:EDE} Assume that for any $T>0$ there exists
a constant $C_{T}>0$ such that for any $|x|\ge1$
\begin{equation}
\label{growth-condition}|a(t,x)|\le C_{T}|x|
\end{equation}
uniformly for $t\in[-T,T]$. Then

\begin{itemize}
\item[(a)] for any initial datum $\psi\in C^{\infty}(\mathbb{R})$ there exists
a unique global (in time) solution $u\in C^{\infty}(\mathbb{R\times\mathbb{R}%
})$;

\item[(b)] if $\psi(x)>0$ $\forall x\in\mathbb{R}$ then $u(t,x)>0$ $\forall
t,x\in\mathbb{R}$.
\end{itemize}
\end{Lemma}

\emph{Proof.} Clearly, the equation (\ref{EDE}) can be rewritten in the form
\begin{equation}
\label{inhomogeneousEq}X(u)=bu
\end{equation}
where $X:=\partial_{t}-a\partial_{x}$. Consider the ordinary differential
equation
\begin{align}
{\dot x}=-a(t,x),\label{1}\\
x|_{t=0}=x_{0}. \label{2}%
\end{align}
It follows from \eqref{growth-condition} that if a solution $x(t,x_{0})$ of
\eqref{1}-\eqref{2} is defined on the interval $t\in(-T,T)$ for some
$0<T<\infty$ then it satisfies the a priory estimate
\[
\sup\limits_{|t|<T}|x(t,x_{0})|<(1+|x_{0}|)\,e ^{C_{T}T}.
\]
In particular, the latter estimate implies that for any $x_{0} \in\mathbb{R},$
there exists a unique global (in time) solution $x(t,x_{0})$ of
\eqref{1}-\eqref{2}. To prove uniqueness of a solution of (\ref{EDE}%
)-(\ref{EDE0}), assume that $u=u(t,x)$ is a smooth solution. It follows from
(\ref{inhomogeneousEq}) that for any $x_{0} \in\mathbb{R}$, the function
$v(t):=u(t,x(t))$ with $x(t):= x(t,x_{0})$ satisfies the differential equation
${\dot v}(t)=b(t,x(t))v(t)$, hence
\begin{equation}
\label{u-evolution}u(t,x(t))=\psi(x_{0})e^{\int_{0}^{t}b(s,x(s))\;ds}.
\end{equation}
As for any given $t \in\mathbb{R},$ the transformation $\mathbb{R}%
\to\mathbb{R}, x_{0}\mapsto x(t,x_{0})$, is a diffeomorphism, formula
(\ref{u-evolution}) defines $u(t,x)$ uniquely. At the same time,
(\ref{u-evolution}) defines a smooth global in time solution of (\ref{EDE}%
)-(\ref{EDE0}). This proves claim $(a)$. Claim $(b)$ also follows from
(\ref{u-evolution}). \hfill$\Box$ \vspace{3mm}

\begin{Rem}
\hspace{-2mm}\textbf{.} Let $a(t,x)=|x|^{\alpha}$ for $|x|\ge1$, where
$\alpha>1$. Then solutions of (\ref{1}) can blow up in finite time. Moreover,
one can show that a solution of \eqref{EDE}-\eqref{EDE0} is not necessarily
unique. This shows that assumption (\ref{growth-condition}) is essential claim
$(a)$ to be true.
\end{Rem}

In the remainder of this appendix, we prove, as advertised, three technical
lemmas used in the proof of Theorem \ref{Th:mKdV-S} and Theorem
\ref{Th:generalized_eigenfunctions}. As above, $\mathcal{S}(\mathbb{R})$
denotes the space of functions $f : \mathbb{R}\to\mathbb{R}$ of Schwartz class.

\begin{Lemma}
\hspace{-2mm}\textbf{.}\label{Lem:EDE-S} Assume that $q(t,x)\in\mathcal{S}%
_{\delta}(\mathbb{R}\times\mathbb{R})$ with $\delta<1$. Then the initial value
problem (\ref{EDE-S1})-(\ref{EDE-S2}) with $\eta\in\mathcal{S}_{-\infty
}(\mathbb{R}\times\mathbb{R})$ and $s_{0}\in\mathcal{S}(\mathbb{R})$ has a
solution in $\mathcal{S}_{-\infty}(\mathbb{R}\times\mathbb{R})$ which is
unique in $C^{\infty}(\mathbb{R}\times\mathbb{R})$.
\end{Lemma}

\emph{Proof.} The initial value problem (\ref{EDE-S1})-(\ref{EDE-S2}) can be
rewritten as
\begin{align}
X(s)  &  =\eta\label{Lie-form1}\\
s|_{t=0}  &  =s_{0} \label{Lie-form2}%
\end{align}
where $X(t,x):=\partial_{t}-2q(t,x)\partial_{x}$ and $X(s)$ denotes the
derivative of $s$ with respect to the flow of the vector field $X$. Denote by
$\xi(t;t_{0},x_{0})$ the solution of the ordinary differential equation
\begin{align}
{\dot x}=-2q(t,x),\label{ODE1}\\
x|_{t=t_{0}}=x_{0} . \label{ODE2}%
\end{align}
If $t_{0}=0$ we denote the corresponding solution $\xi(t;0,x_{0})$ by
$\xi(t,x_{0})$. As $q\in\mathcal{S}_{\delta}(\mathbb{R}\times\mathbb{R})$ and
$\delta<1$ it follows that for any $0<T<\infty$ there exists $C_{T}>0$ such
that for any $|x|\ge1$ and $t\in[-T,T]$
\begin{equation}
\label{growth_condition1}|q(t,x)|\le C_{T}|x|\;.
\end{equation}
In particular, (\ref{growth_condition1}) implies that the solution
$\xi(t;t_{0},x_{0})$ is defined for any $t\in\mathbb{R}$. As $q(t,x)$ is
$C^{\infty}$-smooth in $(t,x)$, the solution $\xi(t;t_{0},x_{0})$ is unique
and depends smoothly on the initial data $(t_{0},x_{0})$. Moreover, for any
given $t_{0},t\in\mathbb{R}$, $t\ge t_{0}$, the transformation $x_{0}%
\mapsto\xi(t;t_{0},x_{0})$, $\mathbb{R}\to\mathbb{R}$, is a diffeomorphism.
Let $s(t,x)$ be a smooth solution of (\ref{Lie-form1})-(\ref{Lie-form2}). Then
the function $s(t):=s(t,\xi(t,x_{0}))$ satisfies the differential equation
${\dot s}=\eta(t,\xi(t,x_{0}))$. In particular,
\begin{equation}
\label{integral_equation}s(t,\xi(t,x_{0}))=s_{0}(x_{0})+\int_{0}^{t}\eta
(\tau,\xi(\tau,x_{0}))\;d\tau.
\end{equation}
Hence, the smooth solution $s(t,x)$ of (\ref{Lie-form1})-(\ref{Lie-form2}) is
defined uniquely by the right side of (\ref{integral_equation}). Equation
(\ref{integral_equation}) can be rewritten in the form
\begin{equation}
\label{integral_equation1}s(t,x)=s_{0}(\xi(0;t,x))+\int_{0}^{t}\eta(\tau
,\xi(\tau;t,x))\;d\tau\;.
\end{equation}
Using that $s_{0}\in\mathcal{S}(\mathbb{R})$, $\eta\in\mathcal{S}_{-\infty
}(\mathbb{R}\times\mathbb{R})$ together with (\ref{integral_equation1}) and
Lemma \ref{Lem:ODE-estimates} $(a)$ stated below one easily gets that for any
$0<T<\infty$ and for any $k\ge0$ there exists a constant $C_{T,k}>0$ such that
for any $t\in[-T,T]$ and any $x$ with $|x|\ge1$
\[
|s(t,x)|\le C_{T,k}|x|^{-k}\;.
\]
Differentiating equation (\ref{integral_equation1}) with respect to $t$ and
$x$, we obtain that for any $k,l\ge0,$ the partial derivative $\partial
^{k}_{t}\partial^{l}_{k}s(t,x)$ is a finite sum
\[
\partial^{k}_{t}\partial^{l}_{k}s(t,x)=\sum_{j}S_{j}(t,x),
\]
where the terms $S_{j}(t,x)$, with the help of Lemma \ref{Lem:ODE-estimates}
below, can be shown to be of the form $S_{j}(t,x)=P_{j}(t,x)Q_{j}(t,x)$ with
$P_{j}\in\mathcal{S}_{-\infty}(\mathbb{R}\times\mathbb{R})$ and $Q_{j}$
growing at most polynomially in $x$ uniformly on compact sets of $t$. In
particular, we get that the solution $s(t,x)$ of the initial value problem
(\ref{EDE-S1})-(\ref{EDE-S2}) lies in $\mathcal{S}_{-\infty}(\mathbb{R}%
\times\mathbb{R})$. The uniqueness of the solution follows from the same
arguments as in in the proof of Lemma \ref{Lem:EDE}. \hfill$\Box$ \vspace{3mm}

The following lemma is used in the proof of Lemma \ref{Lem:EDE-S}. We use the
same notation as in the proof of this lemma.

\begin{Lemma}
\hspace{-2mm}\textbf{.}\label{Lem:ODE-estimates} Assume that $q(t,x)\in
\mathcal{S}_{\delta}(\mathbb{R}\times\mathbb{R})$ with $\delta<1$. Then the
following statements hold:

\begin{itemize}
\item[(a)] For any $0<T<\infty$ there exist constants $C_{1}=C_{1}(T)$,
$C_{2}=C_{2}(T)$, $0<C_{1}<C_{2}$, and $N=N(T)>0$ such that for any
$t,t^{\prime}\in[-T,T]$ and $x$ with $|x|\ge N$
\begin{equation}
\label{estimate1}C_{1}|x|\le|\xi(t;t^{\prime},x)|\le C_{2}|x|\;.
\end{equation}

\item[(b)] For any $0<T<\infty$ and for any $k,l,m\ge0$ with $k+l\ge1,$ there
exists a constant $C_{T,k,l,m}>0$ such that for any $t,t^{\prime}\in[-T,T]$
and $x$ with $|x|\ge1$
\begin{equation}
\label{estimate2}|\partial^{k}_{t}\partial^{l}_{t^{\prime}}\partial^{m}_{x}%
\xi(t;t^{\prime},x)|\le C_{T,k,l,m}|x|^{\delta- m}\;.
\end{equation}

\item[(c)] For any $0<T<\infty$ and for any $m\ge0$ there exists a constant
$C_{T,m}>0$ such that for any $t,t^{\prime}\in[-T,T]$ and $x$ with $|x|\ge1$
\begin{equation}
\label{estimate3}|\partial^{m}_{x}\xi(t;t^{\prime},x)|\le C_{T,m}|x|^{1-m}\;.
\end{equation}

\end{itemize}
\end{Lemma}

\emph{Proof.} Let $R(t,x):=-2q(t,x)$. Clearly,
\begin{equation}
\label{R}R\in\mathcal{S}_{\delta}(\mathbb{R}\times\mathbb{R}),\;\;\delta<1.
\end{equation}

$(a)$ As $R(t,x)$ satisfies for any given $0<T<\infty$ the growth condition
(\ref{growth_condition1}) for $x\ge1$ and $|t|\le T$ with some constant
$C_{T}>0$, the solution $\xi(t;t^{\prime},x)$ (defined globally in time) of
the ordinary differential equation
\begin{eqnarray}
\label{e:ode}{\dot\xi}=R(t,\xi)\label{Eq1}\\
\xi|_{t=t^{\prime}}=x \label{Eq2}%
\end{eqnarray}
satisfies for any $x\ge1$ and $t,t^{\prime}\in[-T,T]$,
\[
-C_{T}\le{\dot\xi}/\xi\le C_{T}%
\]
or
\[
xe^{-C_{T}|t-t^{\prime}|}\le\xi(t;t^{\prime},x)\le xe^{C_{T}|t-t^{\prime}|}.
\]
Hence, for any $x\ge N:=e^{2C_{T}T}$ and $t,t^{\prime}\in[-T,T]$ one has
\[
xe^{-2C_{T}T}\le\xi(t;t^{\prime},x)\le xe^{2C_{T}T}.
\]
Similarly one argues for $x\le-N$ to conclude, altogether, that
\[
e^{-2C_{T}T}|x|\le|\xi(t;t^{\prime},x)|\le e^{2C_{T}T}|x|
\]
for any $t,t^{\prime}\in[-T,T]$ and any $|x|\ge N$.

$(b)$ First define a class of continuous functions $\mathcal{B}^{\delta}%
\equiv\mathcal{B}^{\delta}(\mathbb{R}^{3})$. By definition, a continuous
function $f : \mathbb{R}^{3}\to\mathbb{R}$ is an element in $\mathcal{B}%
^{\delta}(\mathbb{R}^{3})$ iff for any $0<T<\infty$ there exists a constant
$C_{T}>0$ such that for any $t,t^{\prime}\in[-T,T]$ and $|x|\ge1$
\[
|f(t,t^{\prime},x)|\le C_{T}|x|^{\delta}.
\]
We start by proving that for any $k,l,m\ge0$ with $k+l\ge1$ the function
$\partial^{k}_{t}\partial^{l}_{t^{\prime}}\partial^{m}_{x}\xi(t;t^{\prime},x)$
belongs to $\mathcal{B}^{\delta}$. For this purpose it is convenient to
consider instead of (\ref{Eq1})-(\ref{Eq2}) the ordinary differential
equation
\begin{align}
{\dot y}=R(t+t^{\prime},y)\label{Eq11}\\
y|_{t=0}=x \label{Eq22}%
\end{align}
where we consider $t^{\prime}\in\mathbb{R}$ as a parameter. Clearly,
\begin{equation}
\label{equality}y(t;t^{\prime},x)=\xi(t+t^{\prime};t^{\prime},x)
\end{equation}
where $\xi(t;t^{\prime},x)$ is the solution of (\ref{Eq1})-(\ref{Eq2}). Hence
$\partial^{k}_{t}\partial^{l}_{t^{\prime}}\partial^{m}_{x}\xi(t;t^{\prime
},x)\in\mathcal{B}^{\delta}$ if and only if
\begin{equation}
\label{have_to_prove}\partial^{k}_{t}\partial^{l}_{t^{\prime}}\partial^{m}%
_{x}y(t;t^{\prime},x)\in\mathcal{B}^{\delta}.
\end{equation}
We will prove (\ref{have_to_prove}). As, by assumption, $R\in\mathcal{S}%
_{\delta}(\mathbb{R}\times\mathbb{R})$, the equation \eqref{Eq11} together
with the lower and upper bounds in (\ref{estimate1}) imply that $y_{t}%
(t;t^{\prime},x)\in\mathcal{B}^{\delta}$.

Differentiating (\ref{Eq11})-(\ref{Eq22}) with respect to $x$, we obtain that
$y_{x}(t;t^{\prime},x)$ satisfies the differential equation
\begin{align}
(y_{x})_{t}  &  =R_{x}(t+t^{\prime},y)y_{x}\;,\\
y_{x}|_{t=0}  &  =1
\end{align}
hence,
\begin{equation}
\label{rep}y_{x}(t;t^{\prime},x)=e^{\int_{0}^{t}R_{x}(\tau+t^{\prime}%
,y(\tau;t^{\prime},x))\;d\tau}\;.
\end{equation}
As $R_{x}\in\mathcal{S}_{\delta-1}(\mathbb{R}\times\mathbb{R})$ with
$\delta-1<0$ we get from claim (a) that for any $0<T<\infty$ there exists a
constant $C_{T}>0$ such that $\forall t,t^{\prime}\in[-T,T]$ and any
$x\in\mathbb{R}$ one has that
\begin{equation}
\label{y_x}|y_{x}(t;t^{\prime},x)|\le C_{T}\;.
\end{equation}
Analogously, differentiating (\ref{Eq11})-(\ref{Eq22}) with respect to the
variable $t^{\prime}$ one gets
\begin{align}
(y_{t^{\prime}})_{t}  &  =R_{x}(t+t^{\prime},y)y_{t^{\prime}}+R_{t}%
(t+t^{\prime},y)\;,\\
y_{t^{\prime}}|_{t=0}  &  =0\;.
\end{align}
By the method of the variation of parameters, one obtains that
\begin{equation}
\label{y_t'-equation}y_{t^{\prime}}(t;t^{\prime},x)=\Big(\int_{0}^{t}%
b(\tau)e^{-\int_{0}^{\tau}a(u)\;du}\;d\tau\Big)e^{\int_{0}^{t}a(u)\;du}%
\end{equation}
where $a(t)=a(t,t^{\prime},x):=R_{x}(t+t^{\prime},y)$ and $b(t)=b(t,t^{\prime
},x):=R_{t}(t+t^{\prime},y)$. As $R_{x}\in\mathcal{S}_{\delta-1}%
(\mathbb{R}\times\mathbb{R})$ and $R_{t}\in\mathcal{S}_{\delta}(\mathbb{R}%
\times\mathbb{R})$ we get that $a\in\mathcal{B}^{\delta-1}$, $b\in
\mathcal{B}^{\delta}$. Using (\ref{y_t'-equation}) and $\delta-1<0$ one
concludes that $y_{t^{\prime}}\in\mathcal{B}^{\delta}$. Differentiating
successively (\ref{Eq11})-(\ref{Eq22}) with respect to the variables
$t^{\prime}$ and $x$ one obtains an equation of the form
\[
(\partial^{l}_{t^{\prime}}\partial^{m}_{x}y)_{t}=R_{x}(t+t^{\prime
},y)(\partial^{l}_{t^{\prime}}\partial^{m}_{x}y)+B(t,t^{\prime},x)
\]
where the inhomogeneous term $B$ is an element in $\mathcal{B}^{\delta- m}$.
Hence, arguing as above, one concludes that $\partial^{l}_{t^{\prime}}%
\partial^{m}_{x}y(t;t^{\prime},x)\in\mathcal{B}^{\delta- m}$.

In order to prove that
\begin{equation}
\label{t-derivatives}\partial^{k}_{t}(\partial^{l}_{t^{\prime}}\partial
^{m}_{x}y(t;t^{\prime},x))\in\mathcal{B}^{\delta- m}%
\end{equation}
for any $k\ge0$ we use induction in $k$. By the considerations from above
(\ref{t-derivatives}) holds for $k=0$. Assume that $k \ge1$ and $\partial
^{j}_{t}\partial^{l}_{t^{\prime}}\partial^{m}_{x}y\in\mathcal{B}^{\delta- m}$
for $0\le j\le k-1$. Differentiating equation (\ref{Eq11}) with respect to
$t^{\prime}$, $x$, and $t$ we obtain
\begin{align}
(\partial^{k-1}_{t}\partial^{l}_{t^{\prime}}\partial^{m}_{x}y)_{t}  &  =
\partial^{k-1}_{t}\partial^{l}_{t^{\prime}}\partial^{m}_{x}(R(t+t^{\prime
},y))\nonumber\\
&  =R_{x}(t+t^{\prime},y)(\partial^{k-1}_{t}\partial^{l}_{t^{\prime}}%
\partial^{m}_{x}y)+B(t,t^{\prime},x). \label{sum}%
\end{align}
Using (\ref{estimate1}) once again together with the induction hypothesis one
proves that the inhomogeneous term $B$ is in $\mathcal{B}^{\delta- m}$. As
$R_{x}\in\mathcal{B}^{\delta-1}$ with $\delta-1<0$ and $\partial^{k-1}%
_{t}\partial^{l}_{t^{\prime}}\partial^{m}_{x}y\in\mathcal{B}^{\delta- m}$ by
induction hypothesis, formula (\ref{sum}) implies (\ref{t-derivatives}).

$(c)$ Statement $(c)$ follows from \eqref{estimate1} (for $m=0$), \eqref{y_x}
(for $m=1$) and then by differentiating (\ref{rep}) and using that $R_{xx}%
\in\mathcal{S}_{\delta-2}(\mathbb{R}\times\mathbb{R})$ (for $m\ge2$).
\hfill$\Box$ \vspace{3mm}

Let $I\subseteq\mathbb{R}$ be a finite or infinite open interval in
$\mathbb{R}$. Arguing as in the proof of Lemma \ref{Lem:EDE-S}, one proves the
following lemma.

\begin{Lemma}
\hspace{-2mm}\textbf{.}\label{Lem:EDE*} Assume that $a\in\mathcal{O}_{\beta
}(I\times\mathbb{R})$ and $b\in\mathcal{O}_{\beta-1}(I\times\mathbb{R})$ with
$\beta\le1$. Then the initial value problem \eqref{EDE}-\eqref{EDE0} with
$\psi\in\mathcal{S}(\mathbb{R})$ has a solution $u\in\mathcal{S}_{-\infty
}(I\times\mathbb{R})$ that is unique in $C^{\infty}(I\times\mathbb{R})$. In
particular, $u\in C^{1}(I,\mathcal{S}(\mathbb{R}))$.
\end{Lemma}

\section{Appendix: Synopsis of some results on KdV}

For the convenience of the reader, we state in this appendix the results on
the existence and uniqueness of solutions of the KdV equation
\begin{align}
q_{t}-6qq_{x}+q_{xxx}=0\label{KdV1}\\
q|_{t=0}=q_{0} \label{KdV2}%
\end{align}
proved in \cite{B1,BS1,BS2,BS3} which we use in the main body of the paper.

Following earlier work of Menikoff \cite{M}, the authors of \cite{BS1,B1}
prove (among other things) the following theorem:

\begin{Th}
\hspace{-2mm}\textbf{.}\label{Th:SB-as} For any $\beta<1$ and for any initial
data $q_{0}\in\mathcal{S}_{\beta}(\mathbb{R})$ there exists a global (in time)
solution $q\in\mathcal{S}_{\beta}(\mathbb{R}\times\mathbb{R})$ of the initial
value problem (\ref{KdV1})-(\ref{KdV2}).
\end{Th}

Completing results of Menikoff \cite{M} the following uniqueness theorem is
proved in \cite{BS2}, by use of a version of Holmgren's principle.

\begin{Th}
\hspace{-2mm}\textbf{.}\label{Th:SB-uniqueness} For any $T>0,$ there is at
most one solution of (\ref{KdV1})-(\ref{KdV2}) in the classes of functions
$q\in C^{\infty}([0,T]\times\mathbb{R})$ such that
\[
q(t,x)=o(|x|)\;\;\;\mbox{and}\;\;\;\partial_{x}^{k}q(t,x)=O(1)\;\;\;\forall
k\ge1
\]
uniformly in $t\in[0,T]$.
\end{Th}

In \cite{BS3} the following theorems are proved:

\begin{Th}
\hspace{-2mm}\textbf{.}\label{Th:SB-O} For any $\beta<1$ and for any initial
data $q_{0}\in\mathcal{O}_{\beta}(\mathbb{R})$ there exists a global in time
solution $q\in\mathcal{O}_{\beta}(\mathbb{R}\times\mathbb{R})$ of the initial
value problem (\ref{KdV1})-(\ref{KdV2}).
\end{Th}

\begin{Th}
\hspace{-2mm}\textbf{.}\label{Th:SB-o} For any $\beta\le1$ and for any initial
data $q_{0}\in o_{\beta}(\mathbb{R})$ there exists a global in time solution
$q\in o_{\beta}(\mathbb{R}\times\mathbb{R})$ of the initial value problem
(\ref{KdV1})-(\ref{KdV2}).
\end{Th}

\begin{Rem}
\hspace{-2mm}\textbf{.} According to Theorem \ref{Th:SB-uniqueness} the
solutions in Theorem \ref{Th:SB-as}, \ref{Th:SB-O}, and \ref{Th:SB-o} are
unique in the corresponding classes.
\end{Rem}

\begin{Rem}
\hspace{-2mm}\textbf{.} It is likely that the methods developed in
\cite{B1}-\cite{BS3} can be used to prove Theorem \ref{Th:mKdV-S},
\ref{Th:mKdV-O}, and \ref{Th:mKdV-o}. However, the proofs will be much more
difficult than the ones presented in this paper.
\end{Rem}


\begin{thebibliography}{99}                                                                                               %


\bibitem {BeS}F. Berezin, M. Shubin, \emph{The Schr{\"o}dinger Equation},
Kluwer Academic Publishers, 1991

\bibitem {B1}I. Bondareva, \emph{The Korteweg-de Vries equation in classes of
increasing functions with prescribed asymptotics as $|x| \to\infty$\yen },
Math. USSR Sbornik, $\mathbf{50}$(1985), no. 1, 125-135

\bibitem {BS1}I. Bondareva, M. Shubin, \emph{Increasing asymptotic solutions
of the Korteweg-de Vries equation and its higher analogues}, Soviet Math.
Dokl., $\mathbf{26}$(1982), no. 3, 716-719

\bibitem {BS2}I. Bondareva, M. Shubin, \emph{Uniqueness of the solutions of
the Cauchy problem for the Korteweg-de Vries equation in classes of increasing
functions}, Moscow University Mathematics Bulletin, $\mathbf{40}$(1985), no.
3, 53-57

\bibitem {BS3}I. Bondareva, M. Shubin, \emph{Equations of Korteweg - de Vries
type in classes of increasing functions}, Journal of Soviet Mathematics,
$\mathbf{51}$(1990), no. 3, 2323-2332

\bibitem {Dub}B. Dubrovin, \emph{On Hamiltonian perturbations of hyperbolic
systems of conservation laws, II: universality of critical behavior}, to
appear in Comm. in Math. Phys., arXiv:math-ph/0510032 v2

\bibitem {GSS}F. Gesztesy, W. Schweiger, B. Simon, \emph{Commutation methods
applied to the mKdV-equation}, Trans. of AMS, $\mathbf{324}$(2), 1991, 465-525

\bibitem {GS}F. Gesztesy. B. Simon, \emph{Constructing solutions of the mKdV -
equation}, J. Funct. Anal., $\mathbf{89}$(1990), 53-60

\bibitem {GN1}P. G. Grinevich, R. G. Novikov. \emph{Analysis of multisoliton
potentials for the two-dimensional Schr\"{o}dinger equation, and a nonlocal
Riemann problem}, Sov. Math. Dokl. \textbf{33} (1986), 9-12.

\bibitem {GN2}P. G. Grinevich, R. G. Novikov. \emph{Transparent potentials at
fixed energy in dimension two. Fixed-energy dispersion relations for the fast
decaying potentials}, Comm. Math. Phys. \textbf{174} (1995), 409-446.

\bibitem {KPST}T. Kappeler, P. Perry, M. Shubin, P. Topalov, \emph{The Miura
map on the line}, IMRN, $\mathbf{2005:50}$(2005), 3091-3133

\bibitem {KT}T. Kappeler, P. Topalov, \emph{Global well-posedness of mKdV in
$L^{2}(\mathbb{T},\mathbb{R})$}, Commun. in PDE, $\mathbf{30}$(2005), 435-449

\bibitem {KPV}C. Kenig, G. Ponce, L. Vega, \emph{Global solutions for the KdV
equation with unbounded data}, J. Diff. Equ., $\mathbf{139}$(1977), 339-364

\bibitem {Lax}P. Lax, \emph{Integrals of nonlinear equations of evolution and
solitary waves}, Comm. Pure App. Math., $\mathbf{21}$(1968), 467-490

\bibitem {Lax:1974}P. Lax, \emph{Periodic solutions of the KdV equation}, in
\emph{Nonlinear Wave Equations}, ed. A. C. Newall, Proc. of the Summer
Seminar, Potsdam, N.Y., 1972. Lectures in Appl. Math. vol. 15 (1974), 83-96.

\bibitem {Man}S. V. Manakov. \emph{The inverse scattering method and
two-dimensional evolution equations}. Uspekhi Mat. Nauk \textbf{31} (1976),
245-256 (Russian).

\bibitem {Mar}V. Marchenko, \emph{Sturm-Liouville operators and applications},
Operator Theory: Advances and Applications, 22. Birkh\"{a}user Verlag, Basel, 1986

\bibitem {M}A. Menikoff, \emph{The existence of unbounded solutions of the
Korteweg - de Vries equation}, Commun. Pure and Appl. Math., $\mathbf{25}%
$(1972), 407-432

\bibitem {Miura}R. Miura, \emph{Korteweg-de Vries equation and
generalizations. I. A remarkable explicit nonlinear transformation}, J. Math.
Phys., $\mathbf{9}$(1968), 1202-1204

\bibitem {VN}S. P. Novikov, A. P. Veselov. \emph{Finite-zone, two-dimensional
potential Schr\"{o}dinger operators. Explicit formula and evolution
equations}. Sov. Math. Doklady \textbf{30} (1984), 588-591.

\bibitem {Shubin}M. Shubin, \emph{Pseudodifferential operators and spectral
theory}, Springer-Verlag, 1987; Second Edition 2001
\end{thebibliography}
\end{document}